\input ./style/arxiv-general.cfg
\documentclass[MSNbibl,number,citesort,seceqn,dvips]{arxbj}
\makeatletter
   \@ifpackageloaded{graphicx}{}{\usepackage{graphicx}}
\makeatother
\usepackage{multirow,dcolumn}

\aid{0}
\volume{21}
\issue{4}
\pubyear{2015}
\firstpage{2393}
\lastpage{2418}
\doi{10.3150/14-BEJ648} 
\docsubty{FLA}

\makeatletter
\newcommand{\nfrac}[2]{#1/#2}

\newcommand{\vfrac}[2]{(#1)/#2}

\newcolumntype{d}[1]{D{.}{.}{#1}}

\newcommand{\E}{\mathbb{E}}
\newcommand{\PP}{\mathbb{P}}
\newcommand{\N}{\mathbb{N}}

\newcommand{\al}{\alpha}
\newcommand{\be}{\beta}
\newcommand{\si}{\sigma}

\newcommand{\f}{\mathcal F}
\newcommand{\mathcall}{\mathcal L}

\newcommand{\WE}{\widetilde{E}}
\newcommand{\WA}{\widetilde{A}}
\newcommand{\WB}{\widetilde{B}}
\newcommand{\WC}{\widetilde{C}}
\newcommand{\WD}{\widetilde{D}}
\newcommand{\WU}{\widetilde{U}}

\newcommand{\Corr}{\operatorname{Corr}}

\newtheorem{prop}{Proposition}[section]
\newtheorem{cor}[prop]{Corollary}
\newproclaim{defi}[prop]{Definition}
\newtheorem{lem}[prop]{Lemma}
\newtheorem{theo}[prop]{Theorem}
\newremark{rem}[prop]{Remark}
\newproclaim{ass}[prop]{Assumption}
\makeatother

\begin{document}
\begin{frontmatter}

\title{Estimation of integrated volatility of volatility with
applications to goodness-of-fit testing}
\runtitle{Integrated volatility of volatility}

\begin{aug}
\author[A]{\inits{M.}\fnms{Mathias}~\snm{Vetter}\corref{}\ead[label=e1]{mathias.vetter@rub.de}}
\address[A]{Fakult\"{a}t f\"{u}r Mathematik, Ruhr-Universit\"{a}t Bochum, 44780 Bochum,
 Germany.\\ \printead{e1}}
\end{aug}

\received{\smonth{7} \syear{2012}}
\revised{\smonth{3} \syear{2014}}

%
\begin{abstract}
In this paper, we are concerned with nonparametric inference on the
volatility of volatility process in stochastic volatility models. We
construct several estimators for its integrated version in a high-frequency setting, all based on increments of spot volatility
estimators. Some of those are positive by construction, others are bias
corrected in order to attain the optimal rate $n^{-1/4}$. Associated
central limit theorems are proven which can be widely used in practice,
as they are the key to essentially all tools in model validation for
stochastic volatility models. As an illustration we give a brief idea
on a goodness-of-fit test in order to check for a certain parametric
form of volatility of volatility.
\end{abstract}

%
\end{frontmatter}

\section{Introduction} \label{Intro}

Nowadays, stochastic volatility models are standard tools in the
continuous-time modelling of financial time series. Typically, the
underlying (log) price process is assumed to follow a diffusion process
of the form
%
\begin{equation}
\label{darstX} X_t = X_0 + \int_0^t
\mu_s \,\mathrm{d}s + \int_0^t
\sigma_s \,\mathrm{d}W_s,
\end{equation}
where $\mu$ and $\sigma$ can be quite general stochastic processes
themselves. A classical case is where the volatility $\sigma
_s^2=\sigma^2(s,X_s)$ is a function of time and state -- a situation
referred to as the one of a local volatility model. It has turned out
in empirical finance that such models do not fit the data very well, as
some stylised facts such as the leverage effect or volatility
clustering cannot be explained using local volatility only. Stochastic
volatility models, however, are able to reproduce such features, as
they bear an additional source of randomness. In these models, the
volatility process is a diffusion process itself, and we focus on a
rather general situation, namely
%
\begin{equation}
\label{darstsigma} \sigma_t^2 = \sigma_0^2
+ \int_0^t \nu_s \,\mathrm{d}s +
\int_0^t \beta_s \,
\mathrm{d}W_s + \int_0^t
\eta_s \,\mathrm{d}W'_s,
\end{equation}
where $\nu$, $\beta$ and $\eta$ again are suitable stochastic
processes and $W'$ is another Brownian motion, independent of $W$. This
model obviously includes the widely used special case of a volatility
with only one driving Brownian motion, which is $\mathrm{d}\sigma_t^2 = \nu_t
\,\mathrm{d}t + \tau_t \,\mathrm{d}V_t$, where $V$ and $W$ are
 jointly Brownian with some
correlation $\rho$.

Stochastic volatility models are typically parametric ones, and
probably the prime example among those is the Heston model of \cite{Heston}, given by
\begin{eqnarray*}
X_t &=& X_0 + \int_0^t
\biggl(\beta- \frac{\si_s^2}{2} \biggr) \,\mathrm{d}s + \int_0^t
\si_s \,\mathrm{d}W_s,
\\
 \si^2_t
&=& \si^2_0 + \kappa\int_0^t
\bigl(\alpha - \si_s^2\bigr) \,\mathrm{d}s + \xi\int
_0^t \si_s \,\mathrm{d}V_s,
\end{eqnarray*}
for some parameters $\beta, \kappa, \al$ and $\xi$, and with $\Corr
(W,V)=\rho$. Here, the volatility process follows a Cox--Ingersoll--Ross
model, that means it is mean-reverting with mean $\alpha$ and speed
$\kappa$, and both diffusion coefficients are proportional with
parameter $\xi$. Particularly the latter property appears to be rather
typical for stochastic volatility models, and in this sense the Heston
model can be regarded as prototypic. Popular alternatives are for
example coming from the more general (but again parametric) class of
(one factor) CEV models, where the diffusion coefficient $\tau$ of
$\sigma^2$ becomes a general power function of $\sigma$, whereas the
drift part of the volatility remains in principle the same. See \cite{Jones} for a survey.

For this reason, statistical inference for stochastic volatility models
has focused on parametric methods for most times, and usually the
authors provide tools for a specific class of models. However, one is
faced with two severe problems: First, it is in most cases impossible
to assess the distribution of $X$ (or its increments), which makes
standard maximum likelihood theory unavailable. Second, the volatility
process $\sigma^2$ is not observable, and many statistical concepts
have in common that they propose to reproduce the unknown volatility
process from observed option prices, typically by using proxies based
on implied volatility. A survey on early estimation methods in this
context can be found in \cite{CG}. One remarkable exception where
stock price data only is used is the paper of \cite{BZ} who construct
a GMM estimator for the parameters of the Heston model from increments
of realised variance. But also in a general setting with no specific
model in mind, the focus has been on parametric approaches. An early
approach on parameter estimation when $\sigma^2$ is ergodic is the
work of \cite{GC2}, optimal rates are discussed in \cite{Hoffmann}
and \cite{Glot}, and a maximum likelihood approach based on proxies
for the volatility can be found in \cite{AK}. Even nonparametric
concepts have been used to identify parameters of a stochastic
volatility model; see, for example, \cite{BR} or \cite{Vet}.

Genuine nonparametric inference for stochastic volatility models has
typically focused on function estimation. Both \cite{Ren} and \cite{CGR} discuss techniques for the estimation of $f$ and $g$, when the
volatility process satisfies $\mathrm{d}\sigma_t^2 = f(\sigma_t^2)\,\mathrm{d}t +
g(\sigma_t^2)\,\mathrm{d}V_t$. In the more general model-free context of (\ref
{darstsigma}), only \cite{BV} and \cite{WM} have discussed estimation
of functionals of volatility of volatility. While the latter focus on
estimation of a kind of leverage effect which involves the volatility
of volatility process(es), the work of \cite{BV} provides a consistent
estimator for integrated volatility of volatility $\int_0^t \tau_s^2
\,\mathrm{d}s$ in the one-factor case. Their approach is inspired by the
asymptotic behaviour of realised variance, which states that the sum of
squared increments of $\sigma^2$ converges in probability to the
quantity of interest. Since $\sigma^2$ is not observable, the authors
use spot volatility estimators instead.

We will pursue their approach and discuss in detail the asymptotic
behaviour of several estimators for integrated volatility of
volatility, all based on increments of spot volatility estimators, thus
using observations of $X$ only. It turns out that in order to attain
the optimal rate of convergence in this context, it is necessary to
conduct a certain bias correction which destroys positivity of the
estimator -- a feature which is well known from the related problem of
volatility estimation under microstructure noise. Several stable
central limit theorems are provided, and by defining appropriate
estimators for the asymptotic (conditional) variance we obtain feasible
versions as well. The latter results are of theoretical interest on one
hand, but are extremely important from an applied point of view as
well, as they make model validation for stochastic volatility models
possible. Given the tremendous number of such models with entirely
different qualitative behaviours, there is a lack of techniques that
help deciding whether a certain model fits the data appropriately or not.

As a first approach to model validation in this framework, we give a
brief idea on how to do goodness-of-fit testing, but our method is by
no means limited to it. Related procedures can be used to test for
example, whether a Brownian component or jumps are present in the
volatility process and what in general the structure of the jump part
is. Such problems have been solved for the price process $X$ in recent
years (see \cite{JP} for an overview), and in principle the methods
are all based on the estimation of plain integrated volatility $\int_0^t \sigma_s^2 \,\mathrm{d}s$ and further quantities, such as truncated versions
or bipower variation. Using our main results, these concepts can be
translated to the stochastic volatility case by using estimators for
integrated volatility of volatility instead, but usually with the
slower rate of convergence $n^{-1/4}$.

The paper is organised as follows: In Section~\ref{res}, we introduce
our estimators and state the central limit theorems, whereas Section~\ref{mc} is on goodness-of-fit testing in stochastic volatility
models. Some Monte Carlo results can be found in Section~\ref{sim},
followed by some concluding remarks in Section~\ref{conc}. An overview
on some proofs plus a couple of details can be found in the \hyperref[app]{Appendix},
whereas large parts of them have been
relegated to a supplementary article \cite{Vetsup}.

\section{Main results}\label{res}

Let us start with some conditions on the processes involved. All of
these are rather mild and covered by a variety of (stochastic)
volatility models used. The only major restriction is that we will
assume most processes to be continuous for a while and only discuss
briefly later how possible adjustments in order to handle jumps in
price and volatility could look like.

\begin{ass} \label{ass1}
Suppose that the process $X$ is given by (\ref{darstX}), where $W$ is
a standard Brownian motion and the drift process $\mu$ is left
continuous. We assume further that the volatility process $\sigma^2$
is a continuous It\^o semimartingale itself, having the representation
(\ref{darstsigma}). $\nu$ is assumed to be left continuous as well,
whereas $\beta$ satisfies the regularity condition
%
\begin{equation}
\label{darsttau} \beta_s^2 = \beta_0^2
+ \int_0^t \omega_s \,\mathrm{d}s
+ \int_0^t \vartheta ^{(1)}_s
\,\mathrm{d}W_s + \int_0^t
\vartheta^{(2)}_s \,\mathrm{d} W'_s,
\end{equation}
where $\omega$ is locally bounded and each $\vartheta^{(l)}$ is left
continuous, $l=1,2$. A similar condition is assumed to hold for $\eta$
as well. Finally, all processes are defined on the same probability
space $(\Omega,\f,(\f_t)_{t\geq0},\PP)$, and all coefficients are specified in such a way that
$\sigma^2$ is almost surely positive and that $\beta^2$ and $\eta^2$
are either almost surely positive or vanishing identically,
respectively.
\end{ass}

As noted in the \hyperref[Intro]{Introduction}, (\ref{darstsigma}) covers a large class
of volatility models used. For $\eta\equiv0$ we are essentially in
the case of a local volatility model, whereas $\beta_s = \rho\tau_s$
and $\eta_s = \sqrt{1 - \rho^2} \tau_s$ for some process $\tau$
and $\rho\in(-1,1)$ refers to the setting of the typical stochastic
volatility models mentioned before, in which both driving Brownian
motions are correlated with $\rho$. The model in (\ref{darstsigma})
is even more flexible, and it is straight-forward to extend all results
to the case of a multi-factor model driven by more than two independent
Brownian motions as well.

Our aim in the following is to draw inference on the integrated
volatility of volatility up to time $t$, which becomes $\int_0^t
(\beta_s^2 + \eta_s^2) \,\mathrm{d}s$ in our context. Any statistical inference
will be based on high-frequency observations of $X$, and we assume that
the data is recorded at equidistant times. Thus, without loss of
generality let the process be defined on the interval $[0,1]$ and
observed at the time points $i/n$, $i=0, \ldots, n$.

Before we discuss several concepts to assess integrated volatility of
volatility in detail, let us recall the principles of estimation of
standard integrated volatility $\int_0^1 \sigma_s^2 \,\mathrm{d}s$. The usual
estimator in the general model-free setting of (\ref{darstX}) is
realised volatility, given by
\[
RV_t^n = \sum_{i=1}^{\lfloor nt \rfloor}
\bigl(\Delta_i^n X\bigr)^2,
\]
where we set $\Delta_i^n Z = Z_{i/n} - Z_{{(i-1)}/n}$ for any process
$Z$. This estimator is optimal in several respects, even though It\^o
formula proves
%
\begin{equation}
\label{Ito1} \bigl(\Delta_i^n X\bigr)^2 =
\int_{(i-1)/n}^{i/n} \sigma_s^2
\,\mathrm{d}s + 2 \int_{(i-1)/n}^{i/n} (X_s -
X_{(i-1)/n}) \,\mathrm{d}X_s
\end{equation}
only, from which it is simple to see that each squared increment
$(\Delta_i^n X)^2$ is only on average equal to integrated volatility
over the corresponding time interval, but not consistent for it.
(Realised volatility, the sum of the squared increments, however,
\textit{is} consistent for the entire integrated volatility, which is
basically due to a martingale argument.)

Our estimators for integrated volatility of volatility will be based on
a similar intuition: Define statistics via sums of increments such that
each summand is on average equal to integrated volatility of volatility
over the corresponding time interval, but not necessarily consistent.
As before, one would like to build those estimators upon increments of
$\sigma^2$. These are in general not observable, so a proxy for them
is needed. Since we are in a model-free world, a natural estimator for
spot volatility $\sigma_{i/n}^2$ is given by
\[
\label{sigmahat} \hat\sigma_{\nfrac{i}n}^2 = \frac{n}{k_n} \sum
_{j=1}^{k_n} \bigl(\Delta _{i+j}^n
X\bigr)^2,\quad\quad i=0, \ldots, n - k_n,
\]
for some auxiliary (integer-valued) sequence $k_n$. See \cite{APPS} or
\cite{Vet} for details on the asymptotic behaviour of this estimator.
It\^o formula again gives
%
\begin{equation}\label{Ito}
\hat\sigma_{\nfrac{i}n}^2 = \frac{n}{k_n} \sum
_{j=1}^{k_n} 2 \int_{\vfrac{i+j-1}{n}}^{\vfrac{i+j}{n}}(X_s
- X_{\vfrac{i+j-1}{n}}) \,\mathrm{d}X_s + \frac{n}{k_n} \int
_{\nfrac{i}n}^{\vfrac{i+k_n}{n}} \sigma_s^2 \,
\mathrm{d}s =: A_i^n + B_i^n,
\end{equation}
so that $\hat\sigma_{i/n}^2 - \sigma_{i/n}^2$ consists of two
sources of error. From the proofs later on, we see that
$A_i^n=\mathrm{O}_p(\sqrt{1/k_n})$, whereas $B_i^n -\sigma^2_{i/n} = \mathrm{O}_p(\sqrt
{k_n/n})$. Therefore, it appears natural to choose $k_n$ to be of the
order $n^{1/2}$ in order to minimize the error of the spot volatility
estimator (and we will see later that this is indeed the best thing to
do), but we will keep this sequence arbitrary in order to allow for
other estimators as well.

While the choice of the spot volatility estimators depends on the
auxiliary sequence $k_n$, we will introduce a second sequence of
integers $l_n$ which governs the length of the intervals over which
increments of $\hat\sigma^2$ are computed. Thus, the basic element of
our final estimators will be $(\hat\sigma^2_{{(i+l_n)}/{n}} - \hat
\sigma^2_{{i}/{n}})^2$, which can be decomposed as
\[
\bigl(\hat\sigma^2_{{(i+l_n)}/{n}} - \hat\sigma^2_{{i}/{n}}
\bigr)^2 = \bigl(A_{i+l_n}^n - A_i^n
\bigr)^2 + \bigl(B_{i+l_n}^n - B_i^n
\bigr)^2 + 2 \bigl(A_{i+l_n}^n -
A_i^n\bigr) \bigl(B_{i+l_n}^n -
B_i^n\bigr).
\]
The average behaviour of the terms above is discussed in the following
lemma, and it depends crucially on the size of both $k_n$ and $l_n$.

\begin{lem} \label{lem0}
Suppose that Assumption~\ref{ass1} holds and let $\E_i^n[Z]$ denote
conditional expectation of some variable $Z$ with respect to $\mathcal
F_{i/n}$. Set also $M_n=\max(k_n, l_n)$ and $m_n=\min(k_n, l_n)$.
Then we have
\begin{eqnarray*}
\E_i^n\bigl[\bigl(A_{i+l_n}^n -
A_i^n\bigr)^2\bigr] &=& 4 l_n(k_n
M_n)^{-1} \sigma^4_{i/n} \bigl(1 +
\mathrm{O}_p\bigl(M_n^{1/2}n^{-1/2}
\bigr)\bigr),
\\
\E_i^n\bigl[\bigl(B_{i+l_n}^n -
B_i^n\bigr)^2\bigr] &=& l_nm_n(M_n-m_n/3)
(k_nM_nn)^{-1} \bigl(\beta^2_{i/n}
+ \eta^2_{i/n}\bigr) \bigl(1 + \mathrm{O}_p
\bigl(M_n^{1/2}n^{-1/2}\bigr)\bigr).
\end{eqnarray*}
%
\end{lem}

The previous lemma gives us several hints on how to obtain an estimator
for integrated volatility of volatility via sums over $(\hat\sigma
^2_{{(i+l_n)}/{n}} - \hat\sigma^2_{{i}/{n}})^2$. First, information
about $\beta^2_{i/n} + \eta^2_{i/n}$ is contained in increments over
the $B_i^n$ only. Therefore, it appears to be reasonable to choose
$k_n$ and $l_n$ later on in such a way that these terms are at least
not smaller than the bias terms due to increments of~$A_i^n$. Or in
other words, the condition becomes that $n \leq C k_n l_n$ for some
generic $C > 0$.

Also, there are basically two ways to construct an estimator. Either,
pick $k_n$ and $l_n$ such that the bias due to increments of $A_i^n$ is
negligible even after dividing by the rate of convergence. This concept
will lead to the estimator
\[
\hat T_t^n = \sum_{i=0}^{ \lfloor nt  \rfloor- (k_n+l_n)}
k_nM_n\bigl(l_nm_n(M_n-m_n/3)
\bigr)^{-1} \bigl(\hat\sigma^2_{(i+l_n)/n} - \hat
\sigma^2_{i/n}\bigr)^2
\]
which is positive by construction. As noted in the \hyperref[Intro]{Introduction}, this
is the kind of estimator \cite{BV} were looking at. Alternatively, one
can use a bias correction and subtract an estimator for the local
quarticity $\sigma_{i/n}^4$. In this case one loses positivity, but we
will see later that the rate of convergence is much faster in this situation.

Let us pursue the first path for a moment, however. In order to
understand what the rate of convergence for estimation of integrated
volatility of volatility will be, the next result is extremely helpful,
as it gives the central limit theorem for the ``oracle'' estimator
\[
\hat S_t^n = \sum_{i=0}^{ \lfloor nt  \rfloor- (k_n+l_n)}
k_nM_n\bigl(l_nm_n(M_n-m_n/3)
\bigr)^{-1} \bigl(B_{i+l_n}^n - B_i^n
\bigr)^2
\]
which depends on the unobservable increments of $B_i^n$ only. All
results in this section will be pointwise in $t$, even though it is
likely that functional versions hold as well.

\begin{prop} \label{prop1}
Suppose that Assumption~\ref{ass1} holds and that both $k_n \sim c
n^{\al}$ and $l_n \sim d n^{\be}$ hold for some $\al, \be\in(0,1)$
and $c,d > 0$. Let also $M_n$ and $m_n$ be defined as before.
\begin{longlist}[(b)]
\item[(a)] If $\al\neq\be$, we have
\[
\sqrt{\frac{n}{M_n}} \biggl( \hat S_t^n - \int
_0^t \bigl(\beta_s^2 +
\eta _s^2\bigr) \,\mathrm{d}s \biggr) \stackrel{
\mathcall-(s)} {\longrightarrow}\sqrt{4/3} \int_0^t
\bigl(\beta_s^2 + \eta_s^2
\bigr) \,\mathrm{d}\overline W_s.
\]
\item[(b)] For $k_n = l_n$ we have
\[
\sqrt{\frac{n}{M_n}} \biggl( \hat S_t^n - \int
_0^t \bigl(\beta_s^2 +
\eta _s^2\bigr) \,\mathrm{d}s \biggr) \stackrel{
\mathcall-(s)} {\longrightarrow}\sqrt{151/70} \int_0^t
\bigl(\beta_s^2 + \eta_s^2
\bigr) \,\mathrm{d} \overline W_s.
\]
\end{longlist}
In both cases, $\overline W$ is a Brownian motion defined on an
extension of the original probability space and independent of
$\mathcal F$ and the convergence in (\ref{stab1}) is $\mathcal
F$-stable in law.
\end{prop}

\begin{rem}
It is obvious from Proposition~\ref{prop1} that the rate of
convergence becomes faster the smaller $M_n$ is chosen. On the other
hand, the condition $n \leq C k_n l_n$ forces $M_n$ to be at least of
the order $n^{1/2}$. In this case, the rate of convergence in
Proposition~\ref{prop1} becomes $n^{-1/4}$, and this rate is known to
be optimal for this statistical problem. Indeed, a related parametric
setting has been discussed in \cite{Hoffmann} a decade ago, and it was
shown therein that this rate is optimal in the special case, where
$\beta$ vanishes identically and $\eta$ is a function of time and
state, known up to a parameter~$\theta$.
\end{rem}

Our first main theorem specifies conditions for a central limit theorem
for $\hat T_t^n$ and is a simple consequence of Lemma~\ref{lem0} and
Proposition~\ref{prop1}.

\begin{theo} \label{theo2}
Suppose that all the assumptions of Proposition~\ref{prop1} hold true.
If further $n^{3/2} M_n^{-3/2} m_n^{-1} \to0$ and $\alpha\neq\beta
$, then the stable central limit theorem
%
\begin{equation}
\label{clttt} \sqrt{\frac{n}{M_n}} \biggl(\hat T_t^n
- \int_0^t \bigl(\beta_s^2
+ \eta _s^2\bigr) \,\mathrm{d}s \biggr) \stackrel{
\mathcall-(s)} {\longrightarrow}\sqrt{4/3} \int_0^t
\bigl(\beta_s^2 + \eta_s^2
\bigr) \,\mathrm{d} \overline W_s.
\end{equation}
holds true.
\end{theo}

The optimal rate of convergence in this case is obtained for the choice
of $M_n = \mathrm{O}(n^{3/5+\varepsilon})$ and $m_n = \mathrm{O}(n^{3/5})$ and approaches
$n^{-1/5}$ for $\varepsilon\to0$. This proves also that it is no
restriction to assume $\alpha\neq\beta$ above.

In order to obtain an estimator with the optimal rate of convergence,
we choose $l_n$ and $k_n$ to be both the same and of the order
$n^{1/2}$, but as noted above we need a bias correction then.
Therefore, we define with a slight abuse of notation
%
\begin{equation}
\label{esti} \hat R_t^n = \sum
_{i=0}^{ \lfloor nt  \rfloor- 2k_n} \biggl(\frac{3}{2k_n} \bigl(\hat
\sigma_{\vfrac{i+k_n}{n}}^2 - \hat\sigma _{\nfrac{i}{n}}^2
\bigr)^2 - 6 \frac{1}{k_n^2} \hat\sigma_{\nfrac{i}n}^4
\biggr),
\end{equation}
where\vspace*{2pt} $\hat\sigma_{i/n}^4 = \frac{n^2}{3k_n} \sum_{j=1}^{k_n}
|\Delta_{i+j}^n X|^4$ is in general different from $(\hat\sigma
_{i/n}^2)^2$. Its asymptotic behaviour is discussed in the following theorem.

\begin{theo} \label{2.1}
Suppose that Assumption~\ref{ass1} holds and let $k_n = c n^{1/2} +
\mathrm{o}(n^{1/4})$ for some $c>0$. Then
%
\begin{equation}
\label{stab1} \sqrt{\frac{n}{k_n}} \biggl(\hat R_t^n
- \int_0^t \bigl(\beta_s^2
+ \eta _s^2\bigr) \,\mathrm{d}s \biggr) \stackrel{
\mathcall-(s)} {\longrightarrow} U_t
\end{equation}
for all $t > 0$, where the limiting variable has the representation
%
\begin{equation}
\label{Ut} U_t = \int_0^t
\alpha_s \,\mathrm{d} \overline W_s,\quad\quad
\alpha_s^2 = \frac
{48}{c^4} \sigma_s^8
+ \frac{12}{c^2} \sigma_s^4 \bigl(
\beta_s^2 + \eta _s^2\bigr) +
\frac{151}{70} \bigl(\beta_s^2 + \eta_s^2
\bigr)^2.
\end{equation}
\end{theo}

%

\begin{rem}
The situation encountered above has an interesting connection to the
problem of eliminating microstructure noise, as we face similar
problems regarding optimal rates of convergence and positivity of the
estimators. Whereas the optimal rate of convergence for estimating
integrated volatility in the noisy setting is $n^{-1/4}$, standard
estimators attaining this rate are not always positive. To ensure
positivity, one typically accepts a drop in the rate of convergence to
$n^{-1/5}$ as well. See, for example, \cite{BHLS} for a thorough
discussion in a general multivariate setting.
\end{rem}

\begin{rem}
Recently, \cite{JR} discussed efficient estimation of $\int_0^t
g(\sigma^2_s) \,\mathrm{d}s$ for general functions $g$. It turned out that
Riemann sums based on $g(\hat\sigma_{i/n}^2)$ indeed attain the
optimal rate of convergence $n^{-1/2}$ in this context, but again the
choice of $k_n$ affects the limiting distribution. The optimal $k_n
\sim n^{1/2}$ leads to additional bias terms in their setting, and at
least some of these can be avoided by choosing $k_n$ in a different
way.
\end{rem}

The limiting distribution in Theorem~\ref{theo2} and Theorem~\ref
{2.1} is mixed normal, and in order to obtain a feasible central limit
theorem we have to introduce consistent estimators for the respective
conditional variances. These are constructed using the same intuition
as before, and precisely we obtain the following theorem.

\begin{theo} \label{2.4}
 (\textup{a})  Under the conditions of Theorem~\ref{theo2}, we have
\[
\hat Q_t^n = \sum_{i=0}^{ \lfloor nt  \rfloor- (k_n+l_n)}
\frac{4nk_n^2M_n^2}{9(l_nm_n(M_n-m_n/3))^2} \bigl(\hat\sigma ^2_{(i+l_n)/n} - \hat
\sigma^2_{i/n}\bigr)^4 \stackrel{\PP} {
\longrightarrow}\int_0^t \frac{4}3
\bigl(\beta_s^2 + \eta_s^2
\bigr)^2 \,\mathrm{d}s.
\]

 (\textup{b})  In the situation of Theorem~\ref{2.1}, we have
\begin{eqnarray*}
 G^{(1)}_{t,n} &=& \frac{1}n \sum
_{i=1}^{\lfloor nt \rfloor- k_n} \bigl(\hat \sigma_{\nfrac{i}n}^4
\bigr)^2 \stackrel{\PP} {\longrightarrow}\int_0^t
\sigma_s^8 \,\mathrm{d}s,
\\
 G^{(2)}_{t,n} &=& \sum_{i=1}^{\lfloor nt \rfloor-2k_n}
\biggl(\frac
{3}{2k_n} \bigl(\hat\sigma_{\vfrac{i+k_n}{n}}^2 - \hat
\sigma_{\nfrac
{i}{n}}^2\bigr)^2 - 6 \frac{1}{k_n^2} \hat
\sigma_{\nfrac{i}n}^4 \biggr) \hat \sigma_{\nfrac{i}n}^4
\stackrel{\PP} {\longrightarrow}\int_0^t
\sigma_s^4 \bigl(\beta_s^2 +
\eta _s^2\bigr) \,\mathrm{d}s,
\\
 G^{(3)}_{t,n} &=& \sum_{i=1}^{\lfloor nt \rfloor-2k_n}
\frac
{n}{k_n^2} \bigl(\hat\sigma^2_{\vfrac{i+k_n}n} - \hat
\sigma^2_{\nfrac
{i}n}\bigr)^4 \stackrel{\PP} {
\longrightarrow} \int_0^t \biggl(
\frac{48}{c^4}\sigma_s^8 + \frac
{16}{c^2}
\sigma_s^4 \bigl(\beta_s^2 +
\eta_s^2\bigr) + \frac{4}{3} \bigl(\beta
_s^2 + \eta_s^2
\bigr)^2 \biggr) \,\mathrm{d}s.
\end{eqnarray*}
Therefore, as a consequence
\[
\hat P_t^n =\frac{453}{280} G^{(3)}_{t,n}
- \frac{n}{k_n^2} \frac
{486}{35} G^{(2)}_{t,n} -
\frac{n^2}{k_n^4} \frac{1038}{35} G^{(1)}_{t,n} \stackrel{
\PP} {\longrightarrow}\int_0^t
\alpha_s^2 \,\mathrm{d}s.
\]
\end{theo}

\begin{rem}
Theorem~\ref{2.4} shows that a consistent estimator for $\int_0^t
(\beta_s^2 + \eta_s^2)^2 \,\mathrm{d}s$ is for example, given by
\[
\frac{3}{4} G^{(3)}_{t,n} - 12 \frac{n}{k_n^2}
G^{(2)}_{t,n} - 36 \frac{n^2}{k_n^4} G^{(1)}_{t,n},
\]
and its proof suggests that a central limit theorem holds with the same
rate of convergence as before. In general, it is quite likely that this
methods provides estimates for arbitrary even powers of integrated
volatility of volatility. A precise theory is left for future research.
\end{rem}

The properties of stable convergence guarantee that dividing by the
square root of a consistent estimator for the conditional variance
gives a feasible central limit theorem for the estimation of integrated
volatility of volatility. See, for example, \cite{PV} for details.
Therefore, the following corollary can be concluded easily.

\begin{cor} \label{2.5}
 (\textup{a})  Under the assumptions of Theorem~\ref{theo2}, we have for
all $t > 0$
%
\begin{equation}
\label{theo2f} \sqrt{\frac{n}{M_n}} \biggl(\hat T^n_t
- \int_0^t \bigl(\beta_s^2
+ \eta _s^2\bigr) \,\mathrm{d}s \biggr) \bigl(\hat
Q_t^n\bigr)^{-1/2} \stackrel{\mathcall} {
\longrightarrow}\mathcal N(0,1).
\end{equation}

(\textup{b})  Under the assumptions of Theorem~\ref{2.1}, we have for
all $t > 0$
%
\begin{equation}
\label{theo21f} \sqrt{\frac{n}{k_n}} \biggl(\hat R^n_t
- \int_0^t \bigl(\beta_s^2
+ \eta _s^2\bigr) \,\mathrm{d}s \biggr) \bigl(\hat
P_t^n\bigr)^{-1/2} \stackrel{\mathcall} {
\longrightarrow}\mathcal N(0,1).
\end{equation}
\end{cor}

\begin{rem}  So far we have only discussed the case where both
processes have continuous paths. Extensions to the situation of
additional jumps in the price process seem to be possible, but are
already quite involved. The following observation is useful: Whenever
there is a jump within the interval $[i/n,(i+2k_n)/n]$, it appears
squared and blown up by $n^2/k_n^2$ within \mbox{$(\hat\sigma
^2_{{(i+k_n)}/{n}}-\hat\sigma^2_{{i}/{n}})^2$}. This is a much larger
order than the usual $k_n/n$ in the continuous case. For this reason,
it appears as if the truncation method due to \cite{Manc} can be
applied, and a similar intuition holds for the bias correction as well.
Note, however, that the raw statistics in this context are sums of
squared increments of $X$ rather than plain increments of $X$ as for
the power variations encountered in \cite{Manc}. Therefore, the
required techniques are different than the standard ones in this area.

The case of jumps in the volatility appears to be even more
complicated, as these come into play via
\[
B_{i+k_n}^n - B_i^n =
\frac{n}{k_n} \int_{\nfrac{i}n}^{\vfrac
{i+k_n}{n}} \bigl(
\sigma_{s+k_n/n}^2 - \sigma_s^2\bigr) \,
\mathrm{d}s,
\]
and therefore the amount to which each increment is affected by a jump
depends crucially on the time at which the jump occurs. Thus, plain
truncation might not be sufficient in this case and an entirely
different estimator was necessary. Both topics are left for future
research.
\end{rem}

\section{Model checks for stochastic volatility models} \label{mc}

In this section, we propose a first approach to goodness-of-fit testing
for stochastic volatility models. Assume we have representation (\ref
{darstX}) for the log price process $X$, whereas the volatility process
satisfies $\mathrm{d}\sigma_t^2 = \nu_t \,\mathrm{d}t + \tau_t \,\mathrm{d}V_t$
as in typical SV
models. There is still a lot of freedom in the modelling of $\sigma
^2$, and the various proposals in the literature typically differ in
the representation of its diffusion part $\tau$. As noted in the
\hyperref[Intro]{Introduction}, a quite general class of stochastic volatility models is
given by the so-called CEV models, in which $\tau_s^2 = \theta(\sigma
^2_s)^{\gamma}$ for some nonnegative $\gamma$ and an unknown
parameter $\theta$, and the most popular among these is the Heston
model from \cite{Heston}, corresponding to $\gamma= 1$.

In order to construct a test whether a certain functional relationship
between $\sigma$ and $\tau$ is present, we employ a technique which
was already used in \cite{DP} or \cite{VD} when dealing with local
volatility models. Suppose we are interested in testing for $\tau_s^2
= \tau^2(s,X_s,\sigma_s^2,\theta)$, where $\tau^2$ is a given
function and $\theta$ is some unknown (in general multidimensional)
parameter. For simplicity, we will focus on the one-dimensional linear
case only, that is
\[
H_0: \tau_s^2 = \theta\tau^2
\bigl(s,X_s,\sigma_s^2\bigr) \quad\quad\mbox{for all
} s \in[0,1] \mbox{ (a.s.)}
\]
Extensions to the general case follow along the lines of Section~5 in
\cite{VD}.

A test for the null hypothesis will be based on the observation that
$H_0$ is equivalent to $N_t = 0$ for all $t \in[0,1]$ (a.s.), where
the process $N_t$ is given by
\begin{eqnarray*}
N_t &=& \int_0^t \bigl(
\tau_s^2 - \theta_{\mathrm{min}} \tau^2
\bigl(s,X_s,\sigma _s^2\bigr) \bigr) \,
\mathrm{d}s,
\\
 \theta_{\mathrm{min}} &=& \operatorname{arg\,min}\limits
_{\theta} \int
_0^1 \bigl(\tau_s^2 -
\theta\tau^2\bigl(s,X_s,\sigma_s^2
\bigr) \bigr)^2 \,\mathrm{d}s.
\end{eqnarray*}
Assume that the function $\tau^2$ is bounded away from zero. Then a
standard argument from Hilbert space theory shows that $\theta_{\mathrm{min}} =
D^{-1} C$ (and therefore $N_t = R_t - B_t D^{-1} C$), where we have set
$R_t=\int_0^t \tau_s^2 \,\mathrm{d}s$ and
\begin{eqnarray*}
B_t &=& \int_0^t \tau^2
\bigl(s,X_s,\sigma_s^2\bigr) \,\mathrm{d}s,\\
D&=& \int_0^1 \tau ^4
\bigl(s,X_s,\sigma_s^2\bigr) \,\mathrm{d}s,\\
C&=& \int_0^1 \tau_s^2
\tau ^2\bigl(s,X_s,\sigma_s^2
\bigr) \,\mathrm{d}s.
\end{eqnarray*}
To define estimators let $k_n$ as before and recall (\ref{esti}). We set
%
\begin{equation}
\label{taues} \hat\tau_{i/n}^2 = {3n}(2k_n)^{-1}
\bigl(\hat\sigma_{{(i+k_n)}/{n}}^2 - \hat\sigma_{{i}/{n}}^2
\bigr)^2 - 6nk_n^{-2} \hat\sigma_{i/n}^4
\end{equation}
and also $\hat N^n_t = \hat R^n_t - \hat B^n_t (\hat D^n)^{-1} \hat
C^n$ with $\hat R^n_t$ from the previous section, whereas we denote
\begin{eqnarray*}
\hat B^n_t &=& \frac{1}n \sum
_{i=0}^{\lfloor nt \rfloor- k_n} \tau ^2\biggl(
\frac{i}n,X_{\nfrac{i}n},
\hat\sigma_{\nfrac{i}n}^2
\biggr),
\\
\hat D^n&=& \frac{1}n \sum
_{i=0}^{n-k_n} \tau^4\biggl(
\frac{i}n,X_{\nfrac{i}n},\hat\sigma _{\nfrac{i}n}^2
\biggr),
\\
 \hat C^n&=&\frac{1}n \sum
_{i=0}^{n- 2k_n} \hat\tau _{\nfrac{i}n}^2
\tau^2\biggl(\frac{i}n,X_{\nfrac{i}n},\hat
\sigma_{\nfrac{i}n}^2\biggr).
\end{eqnarray*}
In the sequel, we will prove weak convergence of $\hat N^n_t - N_t$, up
to a suitable normalisation. Theorem~\ref{2.1} suggests that $\sqrt
{n/k_n}$ is a reasonable choice, and the following claim proves that
two of the estimators converge at a faster speed, at least if we impose
an additional smoothness condition on the function $\tau^2$.

\begin{lem} \label{lem2}
Suppose that the function $\tau^2$ has continuous partial derivatives
of second order. Then we have
\[
\hat B^n_t - B_t = \mathrm{o}_p
\bigl(n^{-1/4}\bigr), \quad\quad\hat D^n - D = \mathrm{o}_p
\bigl(n^{-1/4}\bigr),
\]
the first result holding uniformly in $t \in[0,1]$.
\end{lem}

The above claim indicates that we have to focus on the terms involving
$\hat\tau_{i/n}^2$ only, which is familiar ground due to the results
of Section~\ref{res}. We start with a proposition on the joint
asymptotic behaviour of $\hat R^n_t$ and $\hat C^n$.

\begin{lem} \label{3.2}
Let $d$ be an integer and $t_1, \ldots, t_d$ be arbitrary in $[0,1]$. Set
\[
\Sigma_{t_1, \ldots, t_d}\bigl(s,X_s,\sigma_s^2
\bigr) = \alpha_s^2 h_{t_1,
\ldots, t_d}
\bigl(s,X_s,\sigma_s^2\bigr)
h_{t_1, \ldots, t_d}\bigl(s,X_s,\sigma_s^2
\bigr)^T
\]
with $h_{t_1, \ldots, t_d}(s,X_s,\sigma_s^2) = (1_{[0,t_1]}, \ldots,
1_{[0,t_d]}, \tau^2(s,X_s,\sigma_s^2))^T$ and $\alpha_s^2$ as in
Theorem~\ref{2.1}. Under the previous assumptions we have the stable
convergence
\[
\sqrt{\frac{n}{k_n}} \bigl(\hat R^n_{t_1} -
R_{t_1}, \ldots, \hat R^n_{t_d} -
R_{t_d}, \hat C^n -C \bigr)^T \stackrel{
\mathcall-(s)} {\longrightarrow} \int_0^1
\Sigma ^{1/2}_{t_1, \ldots, t_d}\bigl(s,X_s,
\sigma_s^2\bigr) \,\mathrm{d}\overline W_s,
\]
where $\overline W$ is a $(d+1)$-dimensional standard Brownian motion
defined on an extension of the original space and independent of
$\mathcal F$.
\end{lem}

We are interested in the asymptotics of the process $A_n(t) = \sqrt
{n/k_n} (\hat N^n_t - N_t)$, and the preceding lemma basically leads to
its finite dimensional convergence. The entire result on weak
convergence of $A_n$ reads as follows.

\begin{theo} \label{3.3}
Assume that the previous assumptions hold. Then the process
$(A_n(t))_{t \in[0,1]}$ converges weakly to a mean zero process
$(A(t))_{t \in[0,1]}$, which is Gaussian conditionally on $\mathcal F$
and whose conditional covariance equals the one of the process
\[
\bigl\{ \alpha_U \bigl(1_{\{U \leq t\}} - B_t
D^{-1} \tau ^2\bigl(U,X_U,
\sigma^2_U\bigr) \bigr) \bigr\}_{t \in[0,1]}
\]
where $U \sim\mathcal U[0,1]$, independent of $\mathcal F$.
\end{theo}

As indicated before, convergence of the finite dimensional
distributions is a direct consequence of Lemma~\ref{3.2}, using the
Delta method for stable convergence (see, e.g., \cite{DPV}). Tightness
follows from Theorem VI. 4.5 in \cite{JS} with a minimal amount of work.

Recall that $N_t=0$ for all $t$ under the null hypothesis. Therefore
Theorem~\ref{3.3} shows that a consistent test is obtained by
rejecting the null hypothesis for large values of a suitable functional
of the process $\{ \sqrt{n/k_n} \hat N^n_t\}_{t \in[0,1]}$. If we
choose the Kolmogorov--Smirnov functional $K_n = \sup_{t \in[0,1]}
\sqrt{n/k_n} |\hat N^n_t|$ for example, we have weak convergence under
the null to $\sup_{t \in[0,1]} |A_t|$ as a consequence of Theorem~\ref{3.3}. The distribution of the latter statistic is extremely
difficult to assess, as it typically depends on the entire process
$(X,\sigma^2)$. We therefore propose to obtain critical values via a
simple bootstrap procedure, which will be introduced in the next section.

To end this section, we define an appropriate estimator for the
conditional variance of $A(t)$, which is given by
\[
s_t^2 = \int_0^t
\alpha_s^2 \,\mathrm{d}s - 2B_tD^{-1}
\int_0^t \alpha_s^2
\tau^2\bigl(s,X_s,\sigma^2_s
\bigr) \,\mathrm{d}s + B_t^2 D^{-2} \int
_0^t \alpha_s^2 \tau
^4\bigl(s,X_s,\sigma^2_s\bigr)
\,\mathrm{d}s,
\]
due to Theorem~\ref{3.3}. Empirical counterparts for $B_t$ and $D$ are
obviously defined by the statistics $\hat B_t$ and $\hat D$, whereas
Theorem~\ref{2.4} suggests that a local estimator for $\alpha
^2_{i/n}$ is given by
\[
\hat\alpha^2_{\nfrac{i}n} = \frac{n^2}{k_n^2} \biggl(
\frac{453}{280} \bigl(\hat\sigma^2_{\vfrac{i+k_n}n} - \hat
\sigma^2_{\nfrac{i}n}\bigr)^4 - \frac{486}{35} \hat
\tau_{\nfrac{i}n}^2 \hat\sigma_{\nfrac{i}n}^4 \biggr)
- \frac{n^6}{k_n^5} \frac{346}{1225} \sum_{j=1}^{k_n}
\bigl|\Delta _{i+j}^n X\bigr|^8.
\]
We obtain the following result, which can be proven in the same way as
Theorem~\ref{2.4}.

\begin{theo} \label{3.4}
Let $t$ be arbitrary and set
\begin{eqnarray*}
\bigl(\hat s^n_t\bigr)^2 &=&
\frac{1}n \sum_{i=1}^{\lfloor nt \rfloor- 2k_n} \hat
\alpha^2_{\nfrac{i}n} - 2\hat B_t \hat D^{-1}
\frac{1}n \sum_{i=1}^{\lfloor nt \rfloor- 2k_n} \hat
\alpha^2_{\nfrac{i}n} \tau^2\biggl(\frac{i}n,X_{\nfrac{i}n},
\hat\sigma^2_{\nfrac{i}n}\biggr)
\\
&&{}+ \hat B_t^2 \hat D^{-2} \frac{1}n
\sum_{i=1}^{\lfloor nt \rfloor- 2k_n} \hat\alpha^2_{\nfrac{i}n}
\tau^2\biggl(\frac{i}n,X_{\nfrac{i}n},\hat\sigma
^2_{\nfrac{i}n}\biggr).
\end{eqnarray*}
Then $(\hat s^n_t)^2$ is consistent for $s_t^2$.
\end{theo}

As a consequence, each statistic $\sqrt{n/k_n} \hat N^n_t/\hat s^n_t$
converges weakly to a normal distribution. This result will be used to
construct a feasible bootstrap statistic in the following.

%

\section{Simulation study} \label{sim}

Let us start with a simulation study concerning the performance of the
rate-optimal $\hat R^n_t$ as an estimator for integrated volatility of
volatility. Throughout this section, we will work with the Heston model
only, and the parameters are chosen as follows: $\beta= 0.3$, $\kappa
= 5$, $\alpha= 0.2$ and $\xi= 0.5$. Furthermore, we set $X_0=0$ and
$\sigma_0^2 = \alpha$. Note that the Feller condition $2\kappa\alpha
\geq\xi^2$ is satisfied, which ensures that the process $\sigma^2$
is almost surely positive as requested. So does $\tau^2$, and it is
obvious that (\ref{darsttau}) holds as well. Therefore all conditions
from Section~\ref{res} are satisfied.

We discuss the finite sample properties of $\hat R^n_t$ for different
choices of the correlation parameter $\rho$ and the number of
observations $n$, and for comparability only we take $n$ to be a square
number and $k_n$ equal to $n^{1/2}$ in all cases, so we have $c=1$.
Theorem~\ref{2.1} suggests that such a medium size of $c$ is
reasonable for finite samples, and additional results not reported here
also point towards the fact that $k_n$ should be chosen close to
$n^{1/2}$. Finally, we set $t=1$. Tables~\ref{tab1}--\ref{tab6}
below are based on 10\,000
simulations.

\begin{table}[t]
\tabcolsep=0pt
 \caption{Mean/variance and simulated quantiles of the feasible
 test statistic (\protect\ref{theo21f}) for $\rho=0$.
 The last column gives the relative amount of negative estimates}\label{tab1}
\begin{tabular*}{\textwidth}{@{\extracolsep{\fill}}ld{2.3}ld{1.4}d{1.4}d{1.4}d{1.4}d{1.4}d{1.4}d{1.4}@{}} \hline
  \multicolumn{1}{@{}l}{$n$} &  \multicolumn{1}{l}{Mean}  &  \multicolumn{1}{l}{Variance} & 0.025 & 0.05 & 0.1 & 0.9 & 0.95 & 0.975 &
  \multicolumn{1}{l@{}}{Neg.}\\
  \hline
  $\hphantom{10\,}400$ & -0.397 & $0.856$ & 0.0497 & 0.0968 & 0.1756 & 0.9754 & 0.9946 & 0.9989 & 0.3522\\
  $\hphantom{0}2\,500$ & -0.287 & $0.965$ & 0.0526 & 0.0932 & 0.1619 & 0.9572 & 0.9862 & 0.9965 & 0.1963\\
  $10\,000$ & -0.170 & $1.023$ & 0.0449 & 0.0799 & 0.1425 & 0.9325 & 0.9757 & 0.9928 & 0.0933 \\
  $22\,500$ & -0.112 & $1.002$  & 0.0404 & 0.0696 & 0.1253  & 0.9271  & 0.9722 & 0.9914 & 0.0510 \\
  $40\,000$ & -0.073 & $1.029$ & 0.0401 & 0.0703 &0.1235 & 0.9203 & 0.9690 & 0.9874 & 0.03 \\
  $52\,900$ & -0.031 & $1.022$ &0.0368 & 0.0653 & 0.1157 & 0.9154 & 0.9633 & 0.9872 & 0.0221 \\ \hline
\end{tabular*}
\end{table}

\begin{table}[b]
\tabcolsep=0pt
\caption{Mean/variance and simulated quantiles of the
feasible test statistic (\protect\ref{theo21f}) for $\rho=-0.2$.
The last column gives the relative amount of negative estimates}\label{tab2}
\begin{tabular*}{\textwidth}{@{\extracolsep{\fill}}ld{2.3}ld{1.4}d{1.4}d{1.4}d{1.4}d{1.4}d{1.4}d{1.4}@{}} \hline
  \multicolumn{1}{@{}l}{$n$} &  \multicolumn{1}{l}{Mean}  &  \multicolumn{1}{l}{Variance} & 0.025 & 0.05 & 0.1 & 0.9 & 0.95 & 0.975 &
  \multicolumn{1}{l@{}}{Neg.}\\
  \hline
  $\hphantom{10\,}400$ & -0.386 & $0.874$ & 0.0491 & 0.0967 & 0.1816 & 0.9724 & 0.9942 & 0.9989 & 0.3528 \\
  $\hphantom{1}2\,500$ & -0.295 & $0.971$ & 0.0552 & 0.0963 & 0.1614 & 0.9559 & 0.9864 & 0.9962 & 0.1996\\
  $10\,000$ & -0.176 & $1.013$ & 0.0464 & 0.0808 & 0.1427& 0.9369 & 0.9770 & 0.9940 & 0.0954\\
  $22\,500$ & -0.226 & $0.987$ & 0.0480 & 0.0840 & 0.1476  & 0.9436  & 0.9776 & 0.9932 & 0.0557\\
  $40\,000$ & -0.075 & $1.001$ & 0.0410 & 0.0673 & 0.1217 & 0.9254 & 0.9713 & 0.9904 & 0.0310 \\
  $52\,900$ & -0.040 & $1.019$ & 0.0396 & 0.0677 & 0.1171 & 0.9180 & 0.9663 & 0.9879 & 0.0246 \\ \hline
\end{tabular*}
\end{table}

Table~\ref{tab1} shows the performance for $\rho= 0$, for which we see that it
takes quite some time for the asymptotics to kick in. Apparent is a
slight overestimation of the lower tails of the distribution, which
seems to originate from the relation of the estimators $\hat R_1$ and
$G^{(3)}_{1,n}$. By construction, in cases where $\hat R_1$ is
underestimating the true quantity, it is typically the case that
increments of $\hat\sigma^2$ are relatively small. As these
increments occur in $G^{(3)}_{1,n}$ as well, most likely the asymptotic
variance is underestimated as well, which explains a too large negative
standardised statistic. The same effect is visible for the upper
quantiles as\vadjust{\goodbreak} well (but resulting in an overestimation), and this simple
explanation is supported by a detailed look at simulation results not
reported here which reveal that the estimation of the asymptotic
variance is extremely accurate for moderate sizes of $\hat R_1 - \int_0^1 \tau_s^2 \,\mathrm{d}s$, but becomes worse when the deviation is rather
large. Similar conclusions can be drawn for the case of a moderately
negative $\rho= -0.2$.

We proceed with the finite sample behaviour of the statistics $\hat
T^n_t$, for which we have a lot of freedom in choosing $k_n$ and $l_n$.
However, in order for both $M_n$ to be rather small and the condition
$n^{3/2} M_n^{-3/2} m_n^{-1} \to0$ to be satisfied, we choose $M_n =
\lfloor n^{3/4} \rfloor$ and $m_n = n^{1/2}$, resulting in a rate of
convergence of about $n^{-1/8}$. Also, we restrict ourselves to $\rho= 0$.

\begin{table}[t]
\tabcolsep=0pt
\caption{Mean/variance and simulated quantiles of the feasible test
statistic (\protect\ref{theo2f})
for $k_n = \lfloor n^{3/4} \rfloor$, $l_n = n^{1/2}$ and $\rho=0$}\label{tab3}
  \begin{tabular*}{\textwidth}{@{\extracolsep{\fill}}ld{2.3}ld{1.4}d{1.4}d{1.4}d{1.4}d{1.4}d{1.4}@{}} \hline
  \multicolumn{1}{@{}l}{$n$} &  \multicolumn{1}{l}{Mean}  &  \multicolumn{1}{l}{Variance} & 0.025 & 0.05 & 0.1 & 0.9 & 0.95 & 0.975 \\
  \hline
  $\hphantom{10\,}400$ & -0.490 & $1.261$ & 0.0966 & 0.1304 & 0.1848 & 0.9990 & 1 & 1 \\
  $\hphantom{1}2\,500$ & -0.320 & $0.837$ & 0.0548 & 0.0837 & 0.1390 & 0.9963 & 0.9999 & 1 \\
  $10\,000$ & -0.291 & $0.806$ & 0.0514 & 0.0800 & 0.1332 & 0.9941 & 1 & 1 \\
  $22\,500$ & -0.259 & $0.920$  & 0.0537 & 0.0873 & 0.1424  & 0.9739  & 0.9948 & 0.9996 \\
  $40\,000$ & -0.215 & $1.040$ & 0.0658 & 0.0920 & 0.1380 & 0.9654 & 0.9970 & 1 \\
  $52\,900$ & -0.164 & $1.016$ & 0.0594 & 0.0826 & 0.1274 & 0.9690 & 0.9988 & 1 \\ \hline
\end{tabular*}
\end{table}

As expected, the approximation of the nominal level is rather poor in
this situation, both when reproducing mean/variance and the quantiles
in the tails. Empirically the results do not improve for other choices
of $k_n$ and $l_n$. Note from Table~\ref{tab3} and Table~\ref{tab4} that results do not
differ very much when choosing either $k_n$ or $l_n$ large, apart from
the remarkable expection of a larger $l_n$ and $n = 10\,000$. But even in
this case, the results are not better than for the rate-optimal $\hat
R^n_t$, which is why we recommend to choose this one rather than $\hat
T^n_t$, even though only the latter estimator is ensured to be positive.

\begin{table}[b]
\tabcolsep=0pt
\caption{Mean/variance and simulated quantiles of the feasible
test statistic (\protect\ref{theo2f})
for $l_n = \lfloor n^{3/4} \rfloor$, $k_n = n^{1/2}$ and $\rho=0$}\label{tab4}
\begin{tabular*}{\textwidth}{@{\extracolsep{\fill}}ld{2.3}ld{1.4}d{1.4}d{1.4}d{1.4}d{1.4}d{1.4}@{}} \hline
  \multicolumn{1}{@{}l}{$n$} &  \multicolumn{1}{l}{Mean}  &  \multicolumn{1}{l}{Variance} & 0.025 & 0.05 & 0.1 & 0.9 & 0.95 & 0.975 \\
  \hline
  $\hphantom{10\,}400$ & -0.476 & $1.255$ & 0.0976 & 0.1316 & 0.1889 & 0.9983 & 0.9999 & 1 \\
  $\hphantom{1}2\,500$ & -0.311 & $0.817$ & 0.0505 & 0.0779 & 0.1322 & 0.9950 & 0.9996 & 1 \\
  $10\,000$ & 0.149 & $1.196$ & 0.0450 & 0.0657 & 0.1005 & 0.8784 & 0.9589 & 0.9904 \\
  $22\,500$ & -0.276 & $0.812$  & 0.0460 & 0.0728 & 0.1234  & 0.9886  & 0.9989 & 1 \\
  $40\,000$ & -0.217 & $1.033$ & 0.0648 & 0.0914 & 0.1354 & 0.9647 & 0.9974 & 1 \\
  $52\,900$ & -0.306 & $0.824$ & 0.0494 & 0.0829 & 0.1456 & 0.9882 & 0.9981 & 1 \\ \hline
\end{tabular*}
\end{table}

As an example for an application in goodness-of-fit testing, we have
constructed a test for a Heston-like volatility structure via a
bootstrap procedure as follows: Based on the\vadjust{\eject} observation that for each
$t$, $\sqrt{n/k_n} \hat N^n_t/\hat s^n_t$ converges weakly to a
standard normal distribution if the null is satisfied, it seems
reasonable to reject the hypothesis for large values of the
standardised Kolmogorov--Smirnov statistic $Y_n = \sup_{i\leq n-2k_n}
|\sqrt{n/k_n} \hat N^n_{i/n}/\hat s^n_{i/n}|$. Since its (asymptotic)
distribution is in general hard to assess, we used bootstrap quantiles
instead, and precisely we have generated bootstrap data
$X^{*(b)}_{i/n}$, $b=1, \ldots, B$, following the equation
\[
X^*_t = \int_0^t
\si^*_s \,\mathrm{d}W^*_s, \quad\quad\bigl(\si_t^*
\bigr)^2 = \hat\alpha+ \int_0^t
\hat\kappa\bigl(\hat\alpha- \bigl(\si_s^*\bigr)^2\bigr)
\,\mathrm{d}s + \hat\xi\int_0^t
\si^*_s \,\mathrm{d}V^*_s.
\]
Here, $W^*$ and $V^*$ are independent Brownian motions, and we have
identified $\hat\alpha$ with the realised volatility of the original
data (which is a measure for the average volatility over $[0,1]$) and
defined $\hat\xi= {\hat\theta}^{1/2}$, since both quantities
coincide under the null. Finally, we have simply set $\hat\kappa= 5
\hat\theta/\hat\alpha$ such that Feller's condition is satisfied.
Setting $B=200$, we have run 500 simulations each.

\begin{table}
\caption{Simulated level of the bootstrap test
based on the standardised Kolmogorov--Smirnov statistic $Y_n$}\label{tab5}
\begin{tabular}{@{}ld{1.3}d{1.3}d{1.3}d{1.3}d{1.3}@{}} \hline
  \multicolumn{1}{@{}l}{$n$} & 0.01 & 0.025 & 0.05 & 0.1 & 0.2\\
  \hline
  $\hphantom{10\,}400$ & 0.004 & 0.012 & 0.024 & 0.064 & 0.172 \\
  $\hphantom{1}2\,500$ & 0.018 & 0.040 & 0.064 & 0.120 & 0.216 \\
  $10\,000$ & 0.010 & 0.018 & 0.040 & 0.084 & 0.194 \\
  $22\,500$ &  0.016 & 0.024 & 0.034 & 0.088 & 0.194 \\
  $40\,000$ & 0.020 & 0.038 & 0.068 & 0.128 & 0.220 \\
  $52\,900$ & 0.010 & 0.020 & 0.052 & 0.118 & 0.200 \\ \hline
\end{tabular}
\end{table}

\begin{table}[b]
\tabcolsep=0pt
\caption{Simulated rejection probabilities of the bootstrap test
based on the standardised
Kolmogorov--Smirnov functional statistic $Y_n$ for various alternatives}\label{tab6}
\begin{tabular*}{\textwidth}{@{\extracolsep{\fill}}ld{1.3}d{1.3}d{1.3}d{1.3}d{1.3}d{1.3}d{1.3}d{1.3}d{1.3}d{1.3}@{}} \hline
\multicolumn{1}{@{}l}{Alt} & \multicolumn{5}{l}{$\gamma=0$} & \multicolumn{5}{l}{$\gamma=2$}  \\[-5pt]
 & \multicolumn{5}{l}{\hrulefill} & \multicolumn{5}{l@{}}{\hrulefill}  \\
\multicolumn{1}{@{}l}{$n$} & 0.01 & 0.025 & 0.05 & 0.1 & 0.2 & 0.01 & 0.025 & 0.05 & 0.1 & 0.2 \\
\hline
$\hphantom{10\,}400$ &  0.032 & 0.072 & 0.124 & 0.192 & 0.292 & 0.056 & 0.080 & 0.128 & 0.204 & 0.320 \\
$\hphantom{1}2\,500$ & 0.028 & 0.052 & 0.082 & 0.134 & 0.262 & 0.044 & 0.090 & 0.156 & 0.248 & 0.372 \\
$10\,000$ & 0.032 & 0.048 & 0.086 & 0.138 & 0.260 & 0.036 & 0.084 & 0.176 & 0.284 & 0.396 \\
$22\,500$ & 0.024 & 0.042 & 0.068 & 0.138 & 0.302 & 0.032 & 0.086 & 0.162 & 0.284 & 0.432 \\
$40\,000$ & 0.028 & 0.046 & 0.094 & 0.196 & 0.426 & 0.028 & 0.064 & 0.120 & 0.310 & 0.482\\
$52\,900$  & 0.026 & 0.040 & 0.082 & 0.174 & 0.422 & 0.024 & 0.058 & 0.144 & 0.320 & 0.488 \\
 \hline
\end{tabular*}
\end{table}

Table~\ref{tab5} shows that the simulated levels are rather close to the
expected ones, irrespectively of~$n$. We have tested two alternatives
from the class of CEV models, namely
\[
\si^2_t = \si^2_0 + \kappa\int
_0^t \bigl(\alpha- \si_s^2
\bigr) \,\mathrm{d}s + V_t
\]
and
\[
 \si^2_t
= \si^2_0 + \kappa\int_0^t
\bigl(\alpha - \si_s^2\bigr) \,\mathrm{d}s + \sqrt{
\kappa} \int_0^t \si_s^2
\,\mathrm{d}V_s,
\]
corresponding to $\gamma=0$ and $\gamma=2$, respectively, and using
the parameters from above. We see from the simulation results that the
rejection probabilities are much larger for the second alternative than
for the first, which can partially explained from two observations:
First, the Vasicek model does not satisfy the assumptions from the
previous sections since the volatility may become negative (in which
case it is set to zero); second, our choice of $\hat\kappa$ is
responsible for a large speed of mean reversion in the bootstrap
algorithm which makes it difficult to distinguish between a Heston-like
volatility of volatility and a constant one. It is expected that the
power improves for an entirely data-driven choice of $\hat\kappa$.

\section{Conclusion} \label{conc}

In this paper, we have discussed a nonparametric method to estimate
the integrated volatility of volatility process in stochastic
volatility models. Our concept is based on spot volatility estimators,
and just as for standard realised volatility we use sums of squares of
these spot volatility estimators to obtain a global estimator for
integrated volatility of volatility. Two classes of estimators have
been investigated -- one consisting of positive estimators with a slow
rate of convergence, the other one being bias corrected but converging
at the optimal rate $n^{-1/4}$. In both cases, central limit theorems
are provided, and we also discuss briefly why a truncated version could
be useful when there are additional jumps in the price process.

Given the variety of stochastic volatility models (in continuous time)
which are used to describe financial data, there is a severe lack in
tools on model validation. Our results fill this gap to a first extent,
as we provide a bootstrap method for goodness-of-fit testing in such
models which investigates whether a specific parametric model for
volatility of volatility is appropriate given the data or not. A
rigorous proof that the proposed procedure keeps the asymptotic level
and is consistent against a large class of alternatives has not been
provided, however, and is left for future research.

%
A different issue to take microstructure issues into account which are
likely to be present when data is observed at high-frequency. Again it
is promising to combine filtering methods for noisy diffusions with the
method proposed in this paper to obtain an estimator for integrated
volatility of volatility in such models as well, but the rate of
convergence is expected to drop further. Precise statements are beyond
the scope of the paper as well.

\begin{appendix}
\section*{Appendix} \label{app}

Note first that every left-continuous process is locally bounded, thus
all processes appearing are. Second, standard localisation procedures
as in \cite{BGJPS} or \cite{J2} allow us to assume that any locally
bounded process is actually bounded, and that almost surely positive
processes can be regarded as bounded away from zero. Universal
constants are denoted by $C$ or $C_r$, the latter if we want to
emphasise dependence on some additional parameter $r$.

Within the main corpus, we give the proof of Theorem~\ref{2.1} only,
which is the by far most complicated result of this work. Analogues of
Lemma~\ref{lem0} and Proposition~\ref{prop1} for the special case of
$l_n = k_n$ are of course parts of it, and it is not difficult to
generalise the proofs in order for both claims to be covered as well.
Therefore, these results are not shown explicitly. Let us start with a
brief sketch of what we will be doing. In general, $\mathcal F$-stable
convergence of a sequence $Z_n$ to some limiting variable $Z$ defined
on an extension $(\widetilde{\Omega},\widetilde{\f},\widetilde{\PP})$ of the original space is equivalent to
%
\begin{equation}
\label{stabc} \E\bigl[h(Z_n) Y\bigr] \to\widetilde\E\bigl[h(Z) Y
\bigr]
\end{equation}
for any bounded Lipschitz function $h$ and any bounded $\mathcal
F$-measurable $Y$. For details, see, for example, \cite{JS} and related
work. Suppose now that there are additional variables $Z_{n,p}$ and
$Z_p$ (the latter defined on the same extension as $Z$) such that
%
\begin{eqnarray}
\label{step1}  \lim_{p\to\infty} \limsup_{n \to\infty}
\E\bigl[|Z_n - Z_{n,p}|\bigr] &=& 0,
\\[-5pt]
\label{step2}  Z_{n,p} &\stackrel{\mathcall-(s)} {\longrightarrow}&
Z_p \quad\quad\mbox{for all } p,
\\
\label{step3}  \lim_{p\to\infty} \widetilde\E\bigl[|Z_p -
Z|\bigr] &=& 0,
\end{eqnarray}
hold. Then the desired stable convergence $Z_n \stackrel{\mathcall-(s)}{\longrightarrow} Z$ follows.
Indeed, let $\varepsilon> 0$. Then there exists a $\delta> 0$ such that
$|x-y| < \delta$ implies $|h(x)-h(y)|< \varepsilon$. Thus we have
\begin{eqnarray*}
&& \bigl| \E\bigl[h(Z_n) Y\bigr] - \E\bigl[h(Z_{n,p}) Y\bigr] \bigr|
\\
&&\quad\leq C \bigl( \E \bigl[h(Z_n)-h(Z_{n,p})|
1_{\{|Z_n - Z_{n,p}| \geq\delta\}}\bigr] + \E \bigl[\bigl|h(Z_n)-h(Z_{n,p})\bigr|
1_{\{|Z_n - Z_{n,p}| < \delta\}}\bigr] \bigr)
\\
&&\quad\leq  C \bigl( P\bigl(|Z_n - Z_{n,p}| \geq\delta\bigr) + \varepsilon
\bigr).
\end{eqnarray*}
We have $\lim_{p\to\infty} \limsup_{n \to\infty} | \E[h(Z_n) Y]
- \E[h(Z_{n,p}) Y] | = 0$ from Markov inequality, (\ref{step1}) and
as $\varepsilon$ was arbitrary. $\lim_{p\to\infty} | \widetilde\E
[h(Z_p) Y] - \widetilde\E[h(Z) Y] | = 0$ can be shown similarly using
(\ref{step3}),
and (\ref{step2}) is by definition equivalent to $
\lim_{n \to\infty} | \E[h(Z_{n,p}) Y] - \widetilde\E[h(Z_p) Y] |
= 0$.
Putting the latter three claims together (plus the triangle inequality
and the fact that all three limiting conditions on $p$ and $n$ are
actually the same) gives (\ref{stabc}).

Our aim in this proof is to employ a certain blocking technique, which
allows us to make use of a type of conditional independence between the
summands within $\hat R_t^n$. To this end, we apply the above
methodology, so we have to define an appropriate double sequence
$U_t^{n,p}$, which will correspond to an approximated version of $\hat
R_t^n$ where we sum over the big blocks only. Some additional notation
is necessary. Let $p \in\N$ be arbitrary. We set
\begin{eqnarray*}
a_{\ell}(p) &=& (\ell-1) (p+2) k_n,
\\
 b_{\ell}(p) &=&
a_\ell(p) + p k_n,
\\
c(p) &=& J_n(p) (p+2)
k_n +1,
\end{eqnarray*}
the first two for any $\ell=1, \ldots, J_n(p)$ with $J_n(p) =
\lfloor{ \lfloor nt-2k_n \rfloor}/{((p+2)k_n)}
\rfloor$. These numbers depend on $n$ as well, even though it does not
show up in the notation. We define further $H_i^n = \int_{(i-1)/{n}}^{i/n} (W_s - W_{{(i-1)}/{n}}) \,\mathrm{d}W_s$. In order to exploit
the afore-mentioned conditional independence, we need approximations
for $A_i^n$ and $B_i^n$ from (\ref{Ito}). For the sake of brevity, we
will only state the approximated increments explicitly, which are given by
%
\begin{eqnarray}
\label{WA}
\nonumber
\WA_{\vfrac{i+k_n}{n}} - \WA_{\nfrac{i}n} &:=&
\frac{n}{k_n} \sum_{j=1}^{k_n} 2
\sigma^2_{\nfrac{a_\ell(p)}{n}} \bigl(H_{i+j+k_n}^n -
H_{i+j}^n\bigr)
\\[-8pt]\\[-8pt]
&=& \frac{n}{k_n} \sigma^2_{\nfrac{a_\ell(p)}{n}} \sum
_{j=1}^{k_n} \bigl( \bigl(\Delta_{i+k_n+j}^n
W\bigr)^2 - \bigl(\Delta_{i+j}^n W
\bigr)^2 \bigr),\nonumber
\end{eqnarray}
where the latter identity is a consequence of It\^o formula, and
%
\begin{eqnarray}
\label{WB} &&\WB_{\vfrac{i+k_n}{n}} - \WB_{\nfrac{i}n}\nonumber
\\[-8pt]\\[-8pt]
&&\quad := \frac{n}{k_n} \int
_{\nfrac{i}n}^{\vfrac{i+k_n}{n}} \bigl(\beta_{\nfrac{a_\ell(p)}{n}}
(W_{s+\nfrac{k_n}{n}}-W_s) + \eta_{\nfrac{a_\ell(p)}{n}} \bigl(W'_{s+\nfrac
{k_n}{n}}-W'_s
\bigr)\bigr) \,\mathrm{d}s.\nonumber
\end{eqnarray}
These quantities are defined for $i=a_\ell(p), \ldots, b_\ell(p)-1$,
thus over the big blocks. For later reasons, we introduce similar
approximations over the small blocks. Set
\begin{eqnarray*}
  \WC_{\vfrac{i+k_n}{n}} - \WC_{\nfrac{i}n} &=& \frac{n}{k_n} \sigma
^2_{\nfrac{b_\ell(p)}{n}} \sum_{j=1}^{k_n}
\bigl( \bigl(\Delta _{i+k_n+j}^n W\bigr)^2 - \bigl(
\Delta_{i+j}^n W\bigr)^2 \bigr),
\\
  \WD_{\vfrac{i+k_n}{n}} - \WD_{\nfrac{i}n} &=& \frac{n}{k_n} \int
_{\nfrac{i}n}^{\vfrac{i+k_n}{n}} \bigl(\beta_{\nfrac{b_\ell(p)}{n}}
(W_{s+\nfrac{k_n}{n}}-W_s) + \eta_{\nfrac{b_\ell(p)}{n}} \bigl(W'_{s+\nfrac
{k_n}{n}}-W'_s
\bigr)\bigr)\,\mathrm{d}s)\,\mathrm{d}s,
\end{eqnarray*}
both for $i=b_\ell(p), \ldots, a_{\ell+1}(p)-1$. Then the following
claim holds, whose proof is postponed to the supplemental file \cite{Vetsup}.

\begin{lem} \label{lem1}
We have
\begin{eqnarray*}
 \E\bigl[ \bigl|A_{\vfrac{i+k_n}{n}} - A_{\nfrac{i}n} - (\WA_{\vfrac{i+k_n}{n}} -
\WA_{\nfrac{i}n})\bigr|^r \bigr]&\leq& C_r \bigl(p
n^{-1}\bigr)^{r/2},
\\
 \E\bigl[\bigl|B_{\vfrac{i+k_n}{n}} - B_{\nfrac{i}n} - (\WB_{\vfrac{i+k_n}{n}} -
\WB_{\nfrac{i}n})\bigr|^r \bigr]&\leq& C_r \bigl(p
n^{-1}\bigr)^{r/2},
\end{eqnarray*}
as well as $\E[|A_{\vfrac{i+k_n}{n}} - A_{\nfrac{i}n}|^r] \leq C_r
n^{-r/4}$ and $\E[|B_{\vfrac{i+k_n}{n}} - B_{\nfrac{i}n}|^r] \leq C_r
n^{-r/4}$ for every $r > 0$. The latter bounds hold also for the
approximated versions, and the same results are true for the
approximation via increments of $\WC$ and $\WD$ over the small blocks.
\end{lem}

Up to a different standardisation, the role of $Z_{n,p}$ in this proof
will be played by $U_t^{n,p} = \sum_{\ell=1}^{J_n(p)} U_\ell^{n,p}$, where
%
\begin{eqnarray}
\label{Uell} U_\ell^{n,p} &=& \sum
_{i=a_\ell(p)}^{b_\ell(p)-1} \frac{3}{2k_n} \bigl((
\WA_{\vfrac{i+k_n}{n}} - \WA_{\nfrac{i}n})+
(\WB_{\vfrac{i+k_n}{n}} - \WB_{\nfrac{i}n})
\bigr)^2 \nonumber
\\[-8pt]\\[-8pt]
&&{}- \frac{pk_n}{n} \biggl[\frac{6n}{k_n^2} \sigma
^4_{\nfrac{a_\ell(p)}{n}} + \bigl(\beta^2_{\nfrac{a_\ell(p)}{n}} + \eta
^2_{\nfrac{a_\ell(p)}{n}}\bigr) \biggr]\nonumber
\end{eqnarray}
involves quantities from the big blocks only. The $U_\ell^{n,p}$ can
be shown to be martingale differences, and the most involved part in
the proof is to use Lemma~\ref{lem1} to obtain
%
\begin{eqnarray}
\label{step1a} \lim_{p\to\infty} \limsup_{n \to\infty}
\sqrt{\frac{n}{k_n}}\E \biggl[ \biggl| \biggl(\hat R^n_t -
\int_0^t \bigl(\beta_s^2
+ \eta_s^2\bigr)\,\mathrm{d}s \biggr) -
U_t^{n,p} \biggr| \biggr]= 0,
\end{eqnarray}
which is the analogue of (\ref{step1}). Let us focus on the remaining
two steps as well. We set
\begin{eqnarray*}
U_t^p &=& \int_0^t
\alpha(p)_s \,\mathrm{d}\overline W_s,
\\
\alpha(p)^2_s &=& \frac{p}{p+2} \biggl(
\frac{48p+d_1}{pc^4} \sigma_s^8 + \frac
{12p+d_2}{pc^2}
\sigma_s^4 \bigl(\beta_s^2 +
\eta_s^2\bigr) + \frac
{151p+d_3}{70p} \bigl(
\beta_s^2 + \eta_s^2
\bigr)^2 \biggr)
\end{eqnarray*}
for certain unspecified constants $d_l$, $l=1,2,3$. In order to prove
the stable convergence
%
\begin{equation}
\label{step2a} \sqrt{\frac{n}{k_n}} U_t^{n,p}
\stackrel{\mathcall-(s)} {\longrightarrow} U_t^p
\end{equation}
we use a well-known result for triangular arrays of martingale
differences, which is due to Jacod \cite{J1}. In particular, the
following three conditions have to be checked.
%
\begin{eqnarray}
\label{step2a1}   \frac{n}{k_n} \sum_{\ell=1}^{J_n(p)}
\E_{a_\ell
(p)}^{n}\bigl[\bigl(U_\ell^{n,p}
\bigr)^2\bigr]& \stackrel{\PP} {\longrightarrow}&\int
_0^t \alpha(p)^2_s \,
\mathrm{d}s,
\\
\label{step2a2}   \frac{n^2}{k_n^2}\sum_{\ell=1}^{J_n(p)}
\E _{a_\ell(p)}^{n}\bigl[\bigl(U_\ell^{n,p}
\bigr)^4\bigr] &\stackrel{\PP} {\longrightarrow}&0,
\\
\label{step2a3}   \sqrt{\frac{n}{k_n}} \sum_{\ell=1}^{J_n(p)}
\E _{a_\ell(p)}^{n}\bigl[U_\ell^{n,p}
(N_{\nfrac{a_{\ell+1}(p)}{n}} - N_{\nfrac{a_{\ell}(p)}{n}})\bigr] &\stackrel{\PP} {\longrightarrow}&0,
\end{eqnarray}
where $N$ is any component of $(W,W')$ or a bounded martingale
orthogonal to both $W$ and $W'$. The final step $
\lim_{p\to\infty} \WE|U_t^p - U_t| = 0$ is obvious.

\subsection{Proof of (\texorpdfstring{\protect\ref{step1a}}{A.8})} For simplicity, we set
$\eta\equiv0$ and $\vartheta^{(2)} \equiv0$ from now on, as
otherwise the proof is exactly the same. In a brief first step, we
replace $\hat R^n_t$ by a version in which the unknown bias and not the
estimator for it is subtracted, that is we introduce
\[
U_t^n = \sum_{i=0}^{ \lfloor nt  \rfloor- 2k_n}
\frac
{3}{2k_n} \bigl(\hat\sigma_{\vfrac{i+k_n}{n}}^2 - \hat
\sigma_{\nfrac
{i}{n}}^2\bigr)^2 - \frac{6}{c^2} \int
_0^t \sigma_s^4 \,
\mathrm{d}s - \int_0^t \beta_s^2
\,\mathrm{d}s.
\]
Theorem~2.1 in \cite{BGJPS} shows that integrals over $\sigma$ can be
estimated with rate $n^{-1/2}$, so the assumption on $k_n$ and a
standard argument regarding boundary terms prove that
\[
\sqrt{\frac{n}{k_n}} \E \biggl[ \biggl| \biggl(\hat R_t^n
- \int_0^t \beta_s^2
\,\mathrm{d}s \biggr) - U_t^n \biggr| \biggr] = \mathrm{o}(1),
\]
uniformly in $t$. A simple consequence of Lemma~\ref{lem1} is that the
remainder terms in $U_t^n$ are negligible, that is
\[
\lim_{p\to\infty} \limsup_{n \to\infty} \sqrt{
\frac{n}{k_n}} \E \Biggl[ \Biggl| \sum_{i=c(p)}^{ \lfloor nt  \rfloor- 2k_n}
\frac{3}{2k_n} \bigl(\hat\sigma_{\vfrac{i+k_n}{n}}^2 - \hat\sigma
_{\nfrac{i}{n}}^2\bigr)^2 - \frac{6}{c^2} \int
_{\nfrac{c(p)}{n}}^t \sigma _s^4 \,
\mathrm{d}s - \int_{\nfrac{c(p)}{n}}^t \beta_s^2
\,\mathrm{d}s \Biggr| \Biggr] = 0,
\]
using also boundedness of the processes on the right hand side and the
definition of $c(p)$. Therefore, we are left to show
%
\begin{equation}
\label{5.1} \lim_{p\to\infty} \limsup_{n \to\infty}
\sqrt{{n}/{k_n}} \E \bigl[ \bigl| \WU_t^{n,p} -
U_t^{n,p} \bigr| \bigr]= 0
\end{equation}
with
%
\begin{eqnarray}
\label{5.1a} \WU_t^{n,p} &=& \sum
_{\ell=1}^{J_n(p)} \Biggl(\sum_{i=a_\ell
(p)}^{b_\ell(p)-1}
+ \sum_{i=b_\ell(p)}^{a_{\ell+1}(p)-1} \Biggr) \bigl(\hat
\sigma_{\vfrac{i+k_n}{n}}^2 - \hat\sigma_{\nfrac{i}{n}}^2
\bigr)^2 \nonumber
\\[-8pt]\\[-8pt]
&&{}- \frac{6}{c^2} \int_0^{\nfrac{c(p)}{n}}
\sigma_s^4 \,\mathrm{d}s - \int_0^{\nfrac{c(p)}{n}}
\beta_s^2\,\mathrm{d}s.\nonumber
\end{eqnarray}

For the integrals within (\ref{5.1a}), recall that these are replaced
by approximated versions in $U_t^{n,p}$. Therefore we have to show for
example,
%
\begin{equation}
\label{appr1} \lim_{p\to\infty} \limsup_{n \to\infty}
\sqrt{\frac{n}{k_n}} \E \Biggl[ \Biggl| \sum_{\ell=1}^{J_n(p)}
\int_{\nfrac{a_\ell
(p)}{n}}^{\nfrac{b_\ell(p)}{n}} \bigl(\beta_s^2
- \beta_{\nfrac{a_\ell
(p)}{n}}^2\bigr)\,\mathrm{d}s \Biggr| \Biggr] = 0.
\end{equation}
For its proof, recall (\ref{darsttau}). The result above follows from
\[
\E \Biggl[ \Biggl| \sum_{\ell=1}^{J_n(p)} \int
_{\nfrac{a_\ell
(p)}{n}}^{\nfrac{b_\ell(p)}{n}} \int_{\nfrac{a_\ell(p)}{n}}^s
\omega _r \,\mathrm{d}r \,\mathrm{d}s \Biggr| \Biggr]\leq C \frac{n}{pk_n}
\biggl(\frac{pk_n}{n}\biggr)^2 \leq C p n^{-1/2}
\]
and
\begin{eqnarray*}
\E \Biggl( \sum_{\ell=1}^{J_n(p)} \int
_{\nfrac{a_\ell(p)}{n}}^{\nfrac
{b_\ell(p)}{n}} \int_{\nfrac{a_\ell(p)}{n}}^s
\vartheta^{(1)}_r \,\mathrm{d}W_r \,
\mathrm{d}s \Biggr)^2 &=& \sum_{\ell=1}^{J_n(p)}
\E \biggl(\int_{\nfrac
{a_\ell(p)}{n}}^{\nfrac{b_\ell(p)}{n}} \int_{\nfrac{a_\ell
(p)}{n}}^s
\vartheta^{(1)}_r \,\mathrm{d}W_r \,
\mathrm{d}s \biggr)^2 
\\
&\leq& C p^2 n^{-1}.
\end{eqnarray*}
Of course, the similar claim
%
\begin{eqnarray}
\label{appr17} \lim_{p\to\infty} \limsup_{n \to\infty}
\sqrt{\frac{n}{k_n}} \E \Biggl[ \Biggl| \sum_{\ell=1}^{J_n(p)}
\int_{\nfrac{a_\ell
(p)}{n}}^{\nfrac{b_\ell(p)}{n}} \bigl(\sigma_s^4
- \sigma_{\nfrac{a_\ell
(p)}{n}}^4\bigr)\,\mathrm{d}s \Biggr| \Biggr] = 0
\end{eqnarray}
holds for the same reasons. We have further
%
\begin{eqnarray}
\label{appr8} \lim_{p\to\infty} \limsup_{n \to\infty}
\sqrt{\frac{n}{k_n}} \E \Biggl[ \Biggl| \sum_{\ell=1}^{J_n(p)}
\frac{pk_n}{n} \biggl(\frac
{n}{k_n^2} - \frac{1}{c^2} \biggr)
\sigma^4_{\nfrac{a_\ell(p)}{n}} \Biggr| \Biggr]= 0,
\end{eqnarray}
which by boundedness of $\sigma$ amounts to prove $n^{-3/4} (k_n^2 -
nc^2) = \mathrm{o}(1)$, and the latter is satisfied by definition of $k_n$. Note
that analogues of (\ref{appr1}), (\ref{appr17}) and (\ref{appr8})
are satisfied over the small blocks as well.

The latter claims prove that we are left to show the approximation over
the big blocks, which is\vspace*{-0.5pt}
%
\begin{eqnarray}
\label{appr3a}
\nonumber
&&\!\!\!\! \lim_{p\to\infty} \limsup
_{n \to\infty} \sqrt{\frac{n}{k_n}} \E \Biggl[ \Biggl| \sum
_{\ell=1}^{J_n(p)} \sum_{i=a_\ell(p)}^{b_\ell
(p)-1}
\frac{3}{2k_n} \bigl(\bigl((A_{\vfrac{i+k_n}{n}} - A_{\nfrac{i}n}) +
(B_{\vfrac{i+k_n}{n}} - B_{\nfrac{i}n})\bigr)^2
\\[-15pt]\\[-2pt]
&&\!\!\!\!\hphantom{\lim_{p\to\infty} \limsup
_{n \to\infty} \sqrt{\frac{n}{k_n}} \E \Biggl[ \Biggl| \sum
_{\ell=1}^{J_n(p)} \sum_{i=a_\ell(p)}^{b_\ell
(p)-1}
\frac{3}{2k_n} \bigl(}{} - \bigl((\WA_{\vfrac{i+k_n}{n}} -
\WA_{\nfrac{i}n})+(\WB_{\vfrac{i+k_n}{n}} -
\WB_{\nfrac{i}n})\bigr)^2 \bigr) \Biggr| \Biggr] = 0,\nonumber
\end{eqnarray}
and the negligibility of the small blocks, that is\vspace*{-0.5pt}
%
\begin{eqnarray}\label{appr7}
&&\lim_{p\to\infty} \limsup_{n \to\infty} \sqrt{
\frac{n}{k_n}} \E \Biggl[ \Biggl| \sum_{\ell=1}^{J_n(p)}
\Biggl(\sum_{i=b_\ell
(p)}^{a_{\ell+1}(p)-1} \frac{3}{2k_n}
\bigl(\hat\sigma_{\vfrac
{i+k_n}{n}}^2 - \hat\sigma_{\nfrac{i}{n}}^2
\bigr)^2\nonumber
\\[-15pt]\\[-2pt]
&&\hphantom{\lim_{p\to\infty} \limsup_{n \to\infty} \sqrt{
\frac{n}{k_n}} \E \Biggl[ \Biggl| \sum_{\ell=1}^{J_n(p)}
\Biggl(}{}- \frac{2k_n}{n} \biggl[\frac{6n}{k_n^2}
\sigma^4_{\vfrac{b_\ell(p)}{n}} + \beta^2_{\nfrac
{b_\ell(p)}{n}} \biggr]
\Biggr) \Biggr| \Biggr] = 0\nonumber
\end{eqnarray}
to obtain (\ref{step1a}).

To prove (\ref{appr3a}), the binomial theorem tells us that we can
discuss the approximation for $B$, the one for $A$ and the mixed part
separately. Using further $x^2 - y^2 = 2y(x-y) + (x-y)^2$ and $xx' -
yy' = (x-y)y' + y(x'-y') + (x-y)(x'-y')$, we see from Lemma~\ref{lem1}
and the growth conditions that (\ref{appr3a}) follows from $\lim_{p\to\infty} \limsup_{n \to\infty} \sum_{r=1}^4 \E[|
L_{n,p}^{(j)} |] = 0$ with\vspace*{-0.5pt}
%
\begin{eqnarray}
 \label{appr3b} L_{n,p}^{(1)} &=& \sqrt{\frac{n}{k_n}} \sum
_{\ell=1}^{J_n(p)} \sum_{i=a_\ell(p)}^{b_\ell(p)-1}
\frac{1}{k_n} \bigl((B_{\vfrac
{i+k_n}{n}} - B_{\nfrac{i}n}) \nonumber
\\[-15pt]\\[-2pt]
&&\hphantom{\sqrt{\frac{n}{k_n}} \sum
_{\ell=1}^{J_n(p)} \sum_{i=a_\ell(p)}^{b_\ell(p)-1}
\frac{1}{k_n} \bigl(}{}- (
\WB_{\vfrac{i+k_n}{n}} - \WB_{\nfrac{i}n}) \bigr) (\WB_{\vfrac{i+k_n}{n}} -
\WB_{\nfrac{i}n}), \nonumber
\\
\label{appr4b}  L_{n,p}^{(2)} &=& \sqrt{\frac{n}{k_n}} \sum
_{\ell=1}^{J_n(p)} \sum_{i=a_\ell(p)}^{b_\ell(p)-1}
\frac{1}{k_n} \bigl((B_{\vfrac
{i+k_n}{n}} - B_{\nfrac{i}n}) \nonumber
\\[-15pt]\\[-2pt]
&&\hphantom{\sqrt{\frac{n}{k_n}} \sum
_{\ell=1}^{J_n(p)} \sum_{i=a_\ell(p)}^{b_\ell(p)-1}
\frac{1}{k_n} \bigl(}{}- (
\WB_{\vfrac{i+k_n}{n}} - \WB_{\nfrac{i}n}) \bigr) (\WA_{\vfrac{i+k_n}{n}} -
\WA_{\nfrac{i}n}),\nonumber
\\
 \label{appr5b} L_{n,p}^{(3)} &=& \sqrt{\frac{n}{k_n}} \sum
_{\ell=1}^{J_n(p)} \sum_{i=a_\ell(p)}^{b_\ell(p)-1}
\frac{1}{k_n} \bigl((A_{\vfrac
{i+k_n}{n}} - A_{\nfrac{i}n}) \nonumber
\\[-15pt]\\[-2pt]
&&\hphantom{\sqrt{\frac{n}{k_n}} \sum
_{\ell=1}^{J_n(p)} \sum_{i=a_\ell(p)}^{b_\ell(p)-1}
\frac{1}{k_n} \bigl(}{}- (
\WA_{\vfrac{i+k_n}{n}} - \WA_{\nfrac{i}n}) \bigr) (\WA_{\vfrac{i+k_n}{n}} -
\WA_{\nfrac{i}n}), \nonumber
\\
\label{appr6b}  L_{n,p}^{(4)} &= &\sqrt{\frac{n}{k_n}} \sum
_{\ell=1}^{J_n(p)} \sum_{i=a_\ell(p)}^{b_\ell(p)-1}
\frac{1}{k_n} \bigl((A_{\vfrac
{i+k_n}{n}} - A_{\nfrac{i}n}) \nonumber
\\[-15pt]\\[-2pt]
&&\hphantom{\sqrt{\frac{n}{k_n}} \sum
_{\ell=1}^{J_n(p)} \sum_{i=a_\ell(p)}^{b_\ell(p)-1}
\frac{1}{k_n} \bigl(}{}- (
\WA_{\vfrac{i+k_n}{n}} - \WA_{\nfrac{i}n}) \bigr) (\WB_{\vfrac{i+k_n}{n}} -
\WB_{\nfrac{i}n}). \nonumber
\end{eqnarray}
Proofs of these claims can be found in the supplementary material \cite{Vetsup}.

Finally, to obtain (\ref{appr7}), we compute the conditional
expectation of the approximated increments, and we will do this for the
$\WA$ and $\WB$ terms only. We have
%
\begin{eqnarray}\label{expa}
&&\E_{a_\ell(p)}^{n}\bigl[(\WA_{\vfrac{i+k_n}{n}} -
\WA_{\nfrac{i}n})^2\bigr]\nonumber
\\
&&\quad = \frac{n^2}{k_n^2} \sigma^4_{\nfrac{a_\ell(p)}{n}}
\sum_{j=1}^{k_n} \E \bigl[ \bigl( \bigl(
\Delta_{i+k_n+j}^n W\bigr)^2 - \bigl(
\Delta_{i+k_n}^n W\bigr)^2 \bigr)^2
\bigr]
\\
&&\quad= \frac{4}{k_n} \sigma^4_{\nfrac{a_\ell(p)}{n}}\nonumber
\end{eqnarray}
as well as
\begin{eqnarray*}
&&\E_{a_\ell(p)}^{n}\bigl[(\WB_{\vfrac{i+k_n}{n}} -
\WB_{\nfrac{i}n})^2\bigr]
\\
&&\quad= 2 \frac{n^2}{k_n^2}
\beta_{\nfrac{a_\ell(p)}{n}}^2 \int_{\nfrac{i}n}^{\vfrac{i+k_n}{n}}
\int_{\nfrac{i}n}^{s} \E\bigl[ (W_{s+\nfrac
{k_n}{n}}-W_s)
(W_{r+\nfrac{k_n}{n}}-W_r)\bigr]\,\mathrm{d}r \,\mathrm{d}s
\\
&&\quad= 2 \frac
{n^2}{k_n^2} \beta_{\nfrac{a_\ell(p)}{n}}^2 \int
_{\nfrac{i}n}^{\vfrac
{i+k_n}{n}} \int_{\nfrac{i}n}^{s}
\biggl(r + \nfrac{k_n}n -s \biggr)\,\mathrm{d}r \,\mathrm{d}s
\\
&&\quad=
\frac{2k_n}{3n} \beta_{\nfrac{a_\ell(p)}{n}}^2.
\end{eqnarray*}
The expectation of the mixed part is zero. Obviously, we have
\begin{eqnarray*}
&&\E_{b_\ell(p)}^{n} \Biggl[ \Biggl(\sum
_{i=b_\ell(p)}^{a_{\ell+1}(p)-1} \frac{3}{2k_n} \bigl((
\WC_{\vfrac{i+k_n}{n}} - \WC_{\nfrac{i}n})+(\WD _{\vfrac{i+k_n}{n}} -
\WD_{\nfrac{i}n})\bigr)^2
\\
&&\hphantom{\E_{b_\ell(p)}^{n} \Biggl[ \Biggl(}{}- \frac{2k_n}{n} \biggl[
\frac
{6n}{k_n^2} \sigma^4_{\nfrac{b_\ell(p)}{n}} + \beta^2_{\nfrac{b_\ell
(p)}{n}}
\biggr] \Biggr) \Biggr] = 0
\end{eqnarray*}
as well. (\ref{appr7}) then follows from the fact that
\begin{eqnarray*}
&&\frac{n}{k_n} \sum_{\ell=1}^{J_n(p)} \E
\Biggl[ \Biggl(\sum_{i=b_\ell
(p)}^{a_{\ell+1}(p)-1}
\frac{3}{2k_n} \bigl((\WC_{\vfrac{i+k_n}{n}} - \WC _{\nfrac{i}n})+(
\WD_{\vfrac{i+k_n}{n}} - \WD_{\nfrac{i}n})\bigr)^2
\\
&&\hphantom{\frac{n}{k_n} \sum_{\ell=1}^{J_n(p)} \E
\Biggl[ \Biggl(}{}-\frac
{2k_n}{n}
\biggl[\frac{6n}{k_n^2} \sigma^4_{\nfrac{b_\ell(p)}{n}} +
\beta^2_{\nfrac{b_\ell(p)}{n}} \biggr] \Biggr)^2 \Biggr]
\end{eqnarray*}
is bounded by a constant times $p^{-1}$, using Lemma~\ref{lem1}.

\subsection{Proof of (\texorpdfstring{\protect\ref{step2a}}{A.9})} Let us check the
conditions for stable convergence in this step, where particularly the
proof of (\ref{step2a1}) is tedious. Write $U_\ell^{n,p} = \sum_{s=1}^3 U_\ell^{n,p,s}$ with
\begin{eqnarray*}
  U_\ell^{n,p,1} &=& \sum_{i=a_\ell(p)}^{b_\ell(p)-1}
\frac{3}{2k_n} \biggl((\WA_{\vfrac{i+k_n}{n}} - \WA_{\nfrac{i}n})^2
- \frac{4}{k_n} \sigma^4_{\nfrac{a_\ell(p)}{n}} \biggr),
\\
  U_\ell^{n,p,2} &=& \sum_{i=a_\ell(p)}^{b_\ell(p)-1}
\frac{3}{2k_n} \biggl((\WB_{\vfrac
{i+k_n}{n}} - \WB_{\nfrac{i}n})^2
- \frac{2k_n}{3n} \beta^2_{\nfrac
{a_\ell(p)}{n}} \biggr),
\\
  U_\ell^{n,p,3} &=& \sum_{i=a_\ell
(p)}^{b_\ell(p)-1}
\frac{3}{k_n} (\WA_{\vfrac{i+k_n}{n}} - \WA _{\nfrac{i}n}) (\WB_{\vfrac{i+k_n}{n}}
- \WB_{\nfrac{i}n}).
\end{eqnarray*}
We have seen in the final step above that these terms are indeed
martingale differences, and it turns out that only the $(U_\ell
^{n,p,s})^2$ terms are responsible for the conditional variance,
whereas the remaining mixed ones are of small order each. To summarize,
the following lemma holds which is proven in the supplementary material
\cite{Vetsup}.

\begin{lem} \label{lem22}
We have
\begin{eqnarray*}
\E_{a_\ell(p)}^{n} \bigl[\bigl(U_\ell^{n,p,1}
\bigr)^2\bigr] &=& \sigma^8_{\nfrac{a_\ell
(p)} {n}} \frac{48p+d_1}{k_n^2} + \mathrm{O}_P\bigl(pn^{-3/2}
\bigr),
\\
\E_{a_\ell(p)}^{n} \bigl[\bigl(U_\ell^{n,p,2}
\bigr)^2\bigr] &=& \beta^4_{\nfrac{a_\ell
(p)} {n}}\biggl(\frac{151}{70}p + d_2\biggr) \frac{k_n^2}{n^2} +
\mathrm{O}_P\bigl(pn^{-3/2}\bigr),
\\
\E_{a_\ell(p)}^{n} \bigl[\bigl(U_\ell^{n,p,3}
\bigr)^2\bigr] &=&\sigma^4_{\nfrac{a_\ell
(p)}{n}}
\beta^2_{\nfrac{a_\ell(p)}{n}} \frac{12p + d_3}{n} + \mathrm{O}_P
\bigl(pn^{-3/2}\bigr),
\end{eqnarray*}
for certain unspecified constants $d_m$, $m=1,2,3$, as well as for each
$r \neq s$
\[
\E_{a_\ell(p)}^{n} \bigl[U_\ell^{n,p,r}
U_\ell^{n,p,s}\bigr] = \mathrm{O}_P
\bigl(pn^{-3/2}\bigr).
\]
\end{lem}

We use Lemma~\ref{lem22} to obtain
\begin{eqnarray*}
&& \frac{n}{k_n} \sum_{\ell=1}^{J_n(p)}
\E_{a_\ell(p)}^{n}\bigl[\bigl(U_\ell ^{n,p}
\bigr)^2\bigr]
\\
&&\quad= \frac{pk_n}n \sum_{\ell=1}^{J_n(p)}
\biggl(\frac
{n^2}{k_n^4} \biggl(48 + \frac{d_1}p\biggr)
\sigma_{\nfrac{a_\ell(p)}{n}}^8
\\
&&\hphantom{\quad= \frac{pk_n}n \sum_{\ell=1}^{J_n(p)}
\biggl(}{}+ \frac{n}{k_n^2} \biggl(12 +
\frac{d_2}{p}\biggr) \sigma_{\nfrac{a_\ell
(p)}{n}}^4
\beta_{\nfrac{a_\ell(p)}{n}}^2 + \biggl(\frac{151}{70} +
\frac
{d_3}p\biggr) \beta_{\nfrac{a_\ell(p)}{n}}^4 \biggr) +
\mathrm{O}_p\biggl(\frac{1}{n^{1/2}}\biggr),
\end{eqnarray*}
thus (\ref{step2a1}) holds using $k_n \sim c n^{1/2}$. Simpler to
obtain is (\ref{step2a2}), as Lemma~\ref{lem1} gives
\[
\frac{n^2}{k_n^2}\sum_{\ell=1}^{J_n(p)}
\E_{a_\ell(p)}^{n}\bigl[\bigl(U_\ell ^{n,p}
\bigr)^4\bigr] \leq C \frac{n^3}{pk_n^3} p^4
n^{-2},
\]
which converges to zero in the usual sense. Finally, one can prove
%
\begin{eqnarray}
\label{finalstep} &&\E_{a_\ell(p)}^{n} \Biggl[ \sum
_{i=a_\ell(p)}^{b_\ell(p)-1} \frac
{3}{2k_n} \bigl((
\WA_{\vfrac{i+k_n}{n}} - \WA_{\nfrac{i}n})\nonumber
\\[-8pt]\\[-8pt]
&&\hphantom{\E_{a_\ell(p)}^{n} \Biggl[ \sum
_{i=a_\ell(p)}^{b_\ell(p)-1} \frac
{3}{2k_n} \bigl(}{}+(\WB_{\vfrac
{i+k_n}{n}} -
\WB_{\nfrac{i}n})\bigr)^2 (N_{\nfrac{a_{\ell+1}(p)}{n}} - N_{\nfrac{a_{\ell}(p)}{n}})
\Biggr] = 0,\nonumber
\end{eqnarray}
where\vspace*{2pt} $N$ is either $W$ or $W'$ or when $N$ is a bounded martingale,
orthogonal to $(W,W')$. Focus on\vspace*{2pt} the first case and decompose $((\WA
_{{(i+k_n)}/{n}} - \WA_{ i/n})+(\WB_{{(i+k_n)}/{n}} - \WB_{i/n}))^2$
via the binomial theorem. For the pure $\WA$ and the pure $\WB$ term,
the claim follows immediately from properties of the normal
distribution upon using that $\sigma_{{a_\ell(p)}/{n}}$ or $\beta
_{{a_\ell(p)}/{n}}$ are $\mathcal F_{{a_\ell(p)}/{n}}$\vspace*{2pt} measurable.
For the mixed term, one has to use the special form\vadjust{\goodbreak} of $\WA
_{{(i+k_n)}/{n}} - \WA_{i/n}$ as a difference of two sums, and a
symmetry argument proves (\ref{finalstep}) in this case. For an
orthogonal $N$, we use standard calculus. By It\^o formula, both $(\WA
_{{(i+k_n)}/{n}} - \WA_{i/n})^2$ and $(\WB_{{(i+k_n)}/{n}} - \WB
_{i/n})^2$ are a measurable variable times the sum of a constant and a
stochastic integral with respect to\vspace*{2pt} $W$ and $W'$, respectively. Thus
(\ref{finalstep}) holds. In the mixed case, we use integration by
parts formula to reduce $(\WA_{{(i+k_n)}/{n}} - \WA_{i/n}) (\WB
_{{(i+k_n)}/{n}} - \WB_{i/n})$ to the sum of a constant, a $\mathrm{d}W$- and a
$\mathrm{d}W'$-integral. Then the same argument applies. Altogether, this gives
(\ref{step2a3}).

\end{appendix}

\section*{Acknowledgements}
The author is grateful for financial support through the collaborative
research center ``Statistical modeling of nonlinear dynamic
processes'' (SFB 823) of the German Research Foundation (DFG). Special
thanks go to two anonymous referees for their valuable comments on
earlier versions of this paper.

\begin{supplement}
\stitle{Additional proofs for claims made in the article}
\slink[doi]{10.3150/14-BEJ648SUPP} 
\sdatatype{.pdf}
\sfilename{BEJ648\_supp.pdf}
\sdescription{We~provide several proofs for either theorems from the
main corpus or additional steps discussed in the \hyperref[app]{Appendix}.}
\end{supplement}

%

\printhistory

\begin{thebibliography}{28}

\bibitem{AK}
%
\begin{barticle}[auto:STB|2014/06/18|12:29:53]
\bauthor{\bsnm{A{\"i}t-Sahalia},~\bfnm{Y.}\binits{Y.}} \AND
\bauthor{\bsnm{Kimmel},~\bfnm{R.}\binits{R.}}
(\byear{2007}).
\btitle{Maximum likelihood estimation of stochastic volatility models}.
\bjournal{J. Financial Economics}
\bvolume{134}
\bpages{507--551}.
\end{barticle}
%
\bptok{imsref}%
\endbibitem

\bibitem{APPS}
%
\begin{barticle}[mr]
\bauthor{\bsnm{Alvarez},~\bfnm{Alexander}\binits{A.}},
\bauthor{\bsnm{Panloup},~\bfnm{Fabien}\binits{F.}},
\bauthor{\bsnm{Pontier},~\bfnm{Monique}\binits{M.}} \AND
\bauthor{\bsnm{Savy},~\bfnm{Nicolas}\binits{N.}}
(\byear{2012}).
\btitle{Estimation of the instantaneous volatility}.
\bjournal{Stat. Inference Stoch. Process.}
\bvolume{15}
\bpages{27--59}.
\bid{doi={10.1007/s11203-011-9062-2}, issn={1387-0874}, mr={2892587}}
\end{barticle}
%
\bptok{imsref}%
\endbibitem

\bibitem{BR}
%
\begin{bmisc}[auto:STB|2014/06/18|12:29:53]
\bauthor{\bsnm{Bandi},~\bfnm{F.}\binits{F.}} \AND
\bauthor{\bsnm{Ren{\`o}},~\bfnm{R.}\binits{R.}}
(\byear{2008}).
\bhowpublished{Nonparametric stochastic volatility.
Technical report}.
\end{bmisc}
%
\bptok{imsref}%
\endbibitem

\bibitem{BV}
%
\begin{bmisc}[auto:STB|2014/06/18|12:29:53]
\bauthor{\bsnm{Barndorff-Nielsen},~\bfnm{O.}\binits{O.}} \AND
\bauthor{\bsnm{Veraart},~\bfnm{A.}\binits{A.}}
(\byear{2009}).
\bhowpublished{Stochastic volatility of volatility in continuous time.
Technical report}.
\end{bmisc}
%
\bptok{imsref}%
\endbibitem

\bibitem{BGJPS}
%
\begin{bincollection}[mr]
\bauthor{\bsnm{Barndorff-Nielsen},~\bfnm{Ole~E.}\binits{O.E.}},
\bauthor{\bsnm{Graversen},~\bfnm{Svend~Erik}\binits{S.E.}},
\bauthor{\bsnm{Jacod},~\bfnm{Jean}\binits{J.}},
\bauthor{\bsnm{Podolskij},~\bfnm{Mark}\binits{M.}} \AND
\bauthor{\bsnm{Shephard},~\bfnm{Neil}\binits{N.}}
(\byear{2006}).
\btitle{A central limit theorem for realised power and bipower
variations of continuous semimartingales}.
In \bbooktitle{From Stochastic Calculus to Mathematical Finance}
\bpages{33--68}.
\blocation{Berlin}:
\bpublisher{Springer}.
\bid{doi={10.1007/978-3-540-30788-4\_3}, mr={2233534}}
\end{bincollection}
%
\bptok{imsref}%
\endbibitem

\bibitem{BHLS}
%
\begin{barticle}[mr]
\bauthor{\bsnm{Barndorff-Nielsen},~\bfnm{Ole~E.}\binits{O.E.}},
\bauthor{\bsnm{Hansen},~\bfnm{Peter~Reinhard}\binits{P.R.}},
\bauthor{\bsnm{Lunde},~\bfnm{Asger}\binits{A.}} \AND
\bauthor{\bsnm{Shephard},~\bfnm{Neil}\binits{N.}}
(\byear{2011}).
\btitle{Multivariate realised kernels: Consistent positive
semi-definite estimators of the covariation of equity prices with noise
and non-synchronous trading}.
\bjournal{J. Econometrics}
\bvolume{162}
\bpages{149--169}.
\bid{doi={10.1016/j.jeconom.2010.07.009}, issn={0304-4076}, mr={2795610}}
\end{barticle}
%
\bptok{imsref}%
\endbibitem

\bibitem{BZ}
%
\begin{barticle}[mr]
\bauthor{\bsnm{Bollerslev},~\bfnm{Tim}\binits{T.}} \AND
\bauthor{\bsnm{Zhou},~\bfnm{Hao}\binits{H.}}
(\byear{2002}).
\btitle{Estimating stochastic volatility diffusion using conditional
moments of integrated volatility}.
\bjournal{J. Econometrics}
\bvolume{109}
\bpages{33--65}.
\bid{doi={10.1016/S0304-4076(01)00141-5}, issn={0304-4076}, mr={1899692}}
\end{barticle}
%
\bptok{imsref}%
\endbibitem

\bibitem{CG}
%
\begin{bincollection}[auto:STB|2014/06/18|12:29:53]
\bauthor{\bsnm{Chernov},~\bfnm{M.}\binits{M.}} \AND
\bauthor{\bsnm{Ghysels},~\bfnm{E.}\binits{E.}}
(\byear{2000}).
\btitle{Estimation of stochastic volatility models for the purpose of
option pricing}.
In \bbooktitle{Computational Finance 1999}
(\beditor{\bfnm{Y.}\binits{Y.}~\bsnm{Abu-Mostafa}},
\beditor{\bfnm{B.}\binits{B.}~\bsnm{LeBaron}},
\beditor{\bfnm{A.}\binits{A.}~\bsnm{Lo}} \AND
\beditor{\bfnm{A.}\binits{A.}~\bsnm{Weigend}}, eds.)
\bpages{567--581}.
\blocation{Cambridge}:
\bpublisher{MIT Press}.
\end{bincollection}
%
\bptok{imsref}%
\endbibitem

\bibitem{CGR}
%
\begin{barticle}[mr]
\bauthor{\bsnm{Comte},~\bfnm{F.}\binits{F.}},
\bauthor{\bsnm{Genon-Catalot},~\bfnm{V.}\binits{V.}} \AND
\bauthor{\bsnm{Rozenholc},~\bfnm{Y.}\binits{Y.}}
(\byear{2010}).
\btitle{Nonparametric estimation for a stochastic volatility model}.
\bjournal{Finance Stoch.}
\bvolume{14}
\bpages{49--80}.
\bid{doi={10.1007/s00780-009-0094-z}, issn={0949-2984}, mr={2563205}}
\end{barticle}
%
\bptok{imsref}%
\endbibitem

\bibitem{DP}
%
\begin{barticle}[mr]
\bauthor{\bsnm{Dette},~\bfnm{Holger}\binits{H.}} \AND
\bauthor{\bsnm{Podolskij},~\bfnm{Mark}\binits{M.}}
(\byear{2008}).
\btitle{Testing the parametric form of the volatility in continuous
time diffusion models -- A stochastic process approach}.
\bjournal{J. Econometrics}
\bvolume{143}
\bpages{56--73}.
\bid{doi={10.1016/j.jeconom.2007.08.002}, issn={0304-4076}, mr={2384433}}
\end{barticle}
%
\bptok{imsref}%
\endbibitem

\bibitem{DPV}
%
\begin{barticle}[mr]
\bauthor{\bsnm{Dette},~\bfnm{Holger}\binits{H.}},
\bauthor{\bsnm{Podolskij},~\bfnm{Mark}\binits{M.}} \AND
\bauthor{\bsnm{Vetter},~\bfnm{Mathias}\binits{M.}}
(\byear{2006}).
\btitle{Estimation of integrated volatility in continuous-time
financial models with applications to goodness-of-fit testing}.
\bjournal{Scand. J. Stat.}
\bvolume{33}
\bpages{259--278}.
\bid{doi={10.1111/j.1467-9469.2006.00479.x}, issn={0303-6898}, mr={2279642}}
\end{barticle}
%
\bptok{imsref}%
\endbibitem

\bibitem{GC2}
%
\begin{barticle}[mr]
\bauthor{\bsnm{Genon-Catalot},~\bfnm{Valentine}\binits{V.}},
\bauthor{\bsnm{Jeantheau},~\bfnm{Thierry}\binits{T.}} \AND
\bauthor{\bsnm{Laredo},~\bfnm{Catherine}\binits{C.}}
(\byear{1999}).
\btitle{Parameter estimation for discretely observed stochastic
volatility models}.
\bjournal{Bernoulli}
\bvolume{5}
\bpages{855--872}.
\bid{doi={10.2307/3318447}, issn={1350-7265}, mr={1715442}}
\end{barticle}
%
\bptok{imsref}%
\endbibitem

\bibitem{Glot}
%
\begin{barticle}[mr]
\bauthor{\bsnm{Gloter},~\bfnm{Arnaud}\binits{A.}}
(\byear{2007}).
\btitle{Efficient estimation of drift parameters in stochastic
volatility models}.
\bjournal{Finance Stoch.}
\bvolume{11}
\bpages{495--519}.
\bid{doi={10.1007/s00780-007-0048-2}, issn={0949-2984}, mr={2335831}}
\end{barticle}
%
\bptok{imsref}%
\endbibitem

\bibitem{Heston}
%
\begin{barticle}[auto:STB|2014/06/18|12:29:53]
\bauthor{\bsnm{Heston},~\bfnm{S.}\binits{S.}}
(\byear{1993}).
\btitle{A closed-form solution for options with stochastic
volatility with applications to bonds and currency options}.
\bjournal{Rev. Financial Studies}
\bvolume{6}
\bpages{327--343}.
\end{barticle}
%
\bptok{imsref}%
\endbibitem

\bibitem{Hoffmann}
%
\begin{barticle}[mr]
\bauthor{\bsnm{Hoffmann},~\bfnm{Marc}\binits{M.}}
(\byear{2002}).
\btitle{Rate of convergence for parametric estimation in a stochastic
volatility model}.
\bjournal{Stochastic Process. Appl.}
\bvolume{97}
\bpages{147--170}.
\bid{doi={10.1016/S0304-4149(01)00130-2}, issn={0304-4149}, mr={1870964}}
\end{barticle}
%
\bptok{imsref}%
\endbibitem

\bibitem{J1}
%
\begin{bincollection}[mr]
\bauthor{\bsnm{Jacod},~\bfnm{Jean}\binits{J.}}
(\byear{1997}).
\btitle{On continuous conditional {G}aussian martingales and stable
convergence in law}.
In \bbooktitle{S\'eminaire de {P}robabilit\'es {XXXI}}.
\bseries{Lecture Notes in Math.}
\bvolume{1655}
\bpages{232--246}.
\blocation{Berlin}:
\bpublisher{Springer}.
\bid{doi={10.1007/BFb0119308}, mr={1478732}}
\end{bincollection}
%
\bptok{imsref}%
\endbibitem

\bibitem{J2}
%
\begin{barticle}[mr]
\bauthor{\bsnm{Jacod},~\bfnm{Jean}\binits{J.}}
(\byear{2008}).
\btitle{Asymptotic properties of realized power variations and related
functionals of semimartingales}.
\bjournal{Stochastic Process. Appl.}
\bvolume{118}
\bpages{517--559}.
\bid{doi={10.1016/j.spa.2007.05.005}, issn={0304-4149}, mr={2394762}}
\end{barticle}
%
\bptok{imsref}%
\endbibitem

\bibitem{JP}
%
\begin{bbook}[mr]
\bauthor{\bsnm{Jacod},~\bfnm{Jean}\binits{J.}} \AND
\bauthor{\bsnm{Protter},~\bfnm{Philip}\binits{P.}}
(\byear{2012}).
\btitle{Discretization of Processes}.
\bseries{Stochastic Modelling and Applied Probability}
\bvolume{67}.
\blocation{Heidelberg}:
\bpublisher{Springer}.
\bid{doi={10.1007/978-3-642-24127-7}, mr={2859096}}
\end{bbook}
%
\bptok{imsref}%
\endbibitem

\bibitem{JR}
%
\begin{barticle}[mr]
\bauthor{\bsnm{Jacod},~\bfnm{Jean}\binits{J.}} \AND
\bauthor{\bsnm{Rosenbaum},~\bfnm{Mathieu}\binits{M.}}
(\byear{2013}).
\btitle{Quarticity and other functionals of volatility: Efficient estimation}.
\bjournal{Ann. Statist.}
\bvolume{41}
\bpages{1462--1484}.
\bid{doi={10.1214/13-AOS1115}, issn={0090-5364}, mr={3113818}}
\end{barticle}
%
\bptok{imsref}%
\endbibitem

\bibitem{JS}
%
\begin{bbook}[mr]
\bauthor{\bsnm{Jacod},~\bfnm{Jean}\binits{J.}} \AND
\bauthor{\bsnm{Shiryaev},~\bfnm{Albert~N.}\binits{A.N.}}
(\byear{2003}).
\btitle{Limit Theorems for Stochastic Processes},
\bedition{2nd} ed.
\bseries{Grundlehren der Mathematischen Wissenschaften [Fundamental
Principles of Mathematical Sciences]}
\bvolume{288}.
\blocation{Berlin}:
\bpublisher{Springer}.
\bid{doi={10.1007/978-3-662-05265-5}, mr={1943877}}
\end{bbook}
%
\bptok{imsref}%
\endbibitem

\bibitem{Jones}
%
\begin{barticle}[mr]
\bauthor{\bsnm{Jones},~\bfnm{Christopher~S.}\binits{C.S.}}
(\byear{2003}).
\btitle{The dynamics of stochastic volatility: Evidence from
underlying and options markets}.
\bjournal{J. Econometrics}
\bvolume{116}
\bpages{181--224}.
\bnote{Frontiers of financial econometrics and financial engineering}.
\bid{doi={10.1016/S0304-4076(03)00107-6}, issn={0304-4076}, mr={2002525}}
\end{barticle}
%
\bptok{imsref}%
\endbibitem

\bibitem{Manc}
%
\begin{barticle}[mr]
\bauthor{\bsnm{Mancini},~\bfnm{Cecilia}\binits{C.}}
(\byear{2009}).
\btitle{Non-parametric threshold estimation for models with stochastic
diffusion coefficient and jumps}.
\bjournal{Scand. J. Stat.}
\bvolume{36}
\bpages{270--296}.
\bid{doi={10.1111/j.1467-9469.2008.00622.x}, issn={0303-6898}, mr={2528985}}
\end{barticle}
%
\bptok{imsref}%
\endbibitem

\bibitem{PV}
%
\begin{barticle}[mr]
\bauthor{\bsnm{Podolskij},~\bfnm{Mark}\binits{M.}} \AND
\bauthor{\bsnm{Vetter},~\bfnm{Mathias}\binits{M.}}
(\byear{2010}).
\btitle{Understanding limit theorems for semimartingales: A short survey}.
\bjournal{Stat. Neerl.}
\bvolume{64}
\bpages{329--351}.
\bid{doi={10.1111/j.1467-9574.2010.00460.x}, issn={0039-0402}, mr={2683464}}
\end{barticle}
%
\bptok{imsref}%
\endbibitem

\bibitem{Ren}
%
\begin{barticle}[mr]
\bauthor{\bsnm{Ren{\`o}},~\bfnm{Roberto}\binits{R.}}
(\byear{2006}).
\btitle{Nonparametric estimation of stochastic volatility models}.
\bjournal{Econom. Lett.}
\bvolume{90}
\bpages{390--395}.
\bid{doi={10.1016/j.econlet.2005.09.009}, issn={0165-1765}, mr={2212176}}
\end{barticle}
%
\bptok{imsref}%
\endbibitem

\bibitem{Vet}
%
\begin{barticle}[mr]
\bauthor{\bsnm{Vetter},~\bfnm{Mathias}\binits{M.}}
(\byear{2012}).
\btitle{Estimation of correlation for continuous semimartingales}.
\bjournal{Scand. J. Stat.}
\bvolume{39}
\bpages{757--771}.
\bid{doi={10.1111/j.1467-9469.2012.00783.x}, issn={0303-6898}, mr={3000847}}
\bptnote{check year}%
\end{barticle}
%
\bptok{imsref}%
\endbibitem

\bibitem{Vetsup}
%
\begin{bmisc}[auto:STB|2014/06/18|12:29:53]
\bauthor{\bsnm{Vetter},~\bfnm{M.}\binits{M.}}
(\byear{2014}).
\bhowpublished{Supplement to ``Estimation of integrated
volatility of volatility with applications to goodness-of-fit testing.''
DOI:\doiurl{10.3150/14-BEJ648SUPP}}.
\end{bmisc}
%
\bptok{imsref}%
\endbibitem

\bibitem{VD}
%
\begin{barticle}[mr]
\bauthor{\bsnm{Vetter},~\bfnm{Mathias}\binits{M.}} \AND
\bauthor{\bsnm{Dette},~\bfnm{Holger}\binits{H.}}
(\byear{2012}).
\btitle{Model checks for the volatility under microstructure noise}.
\bjournal{Bernoulli}
\bvolume{18}
\bpages{1421--1447}.
\bid{doi={10.3150/11-BEJ384}, issn={1350-7265}, mr={2995803}}
\bptnote{check year}%
\end{barticle}
%
\bptok{imsref}%
\endbibitem

\bibitem{WM}
%
\begin{barticle}[mr]
\bauthor{\bsnm{Wang},~\bfnm{Christina~D.}\binits{C.D.}} \AND
\bauthor{\bsnm{Mykland},~\bfnm{Per~A.}\binits{P.A.}}
(\byear{2014}).
\btitle{The {e}stimation of {l}everage {e}ffect {w}ith
{h}igh-{f}requency {d}ata}.
\bjournal{J.~Amer. Statist. Assoc.}
\bvolume{109}
\bpages{197--215}.
\bid{doi={10.1080/01621459.2013.864189}, issn={0162-1459}, mr={3180557}}
\end{barticle}
%
\bptok{imsref}%
\endbibitem

\end{thebibliography}
\end{document}